\DeclareFontFamily{U}{wncy}{}
\DeclareFontShape{U}{wncy}{m}{n}{%
   <5>wncyr5%
   <6>wncyr6%
   <7>wnyr7%
   <8>wncyr8%
   <9>wncyr9%
   <10>wncyr10%
   <11>wncyr10%
   <12>wncyr6%
   <14>wncyr7%
   <17>wncyr8%
   <20>wncyr10%
   <25>wncyr10}{}

\documentclass[a4paper,11pt]{article}

\usepackage{latexsym,graphics,color}
\usepackage{amsmath,amsfonts,amssymb,amsthm}
\newtheorem{thm}{Theorem}[section]
\newtheorem{lem}[thm]{Lemma}

\newtheorem{cor}[thm]{Corollary}
\newtheorem{prop}[thm]{Proposition}

\newtheorem{rem}[thm]{Remark}

\input epsf

\title{Constructions of some perfect integral lattices with 
minimum $4$}
\author{Roland Bacher\footnote{This work has been partially supported by the LabEx PERSYVAL-Lab (ANR--11-LABX-0025). The author is a member of the project-team GALOIS supported by this LabEx.}}

\begin{document}
\maketitle

{\sl Abstract\footnote{Keywords: Perfect lattice, finite abelian group,
projective plane,  
equiangular system, Schl\"afli graph, Sidon set, Craig lattice.
Math. class:  Primary: 11H55, Secondary: 11T06, 20K01, 05B30, 05E30.}: 
We construct several families of perfect 
sublattices with minimum $4$ of $\mathbb Z^d$. In particular,
the number of $d$-dimensional perfect integral lattices with minimum $4$
grows faster than $d^k$ for every exponent $k$.
}
\vskip.5cm

\section{Perfection and perfect lattices}

A subset $\mathcal S$ of a real $d$-dimensional vector space 
$V$ is a \emph{perfect subset of $V$} 
(or perfect in $V$)
if the span of the set 
$\{v\otimes v\}_{v\in \mathcal S}$ is the full 
${d+1\choose 2}$-dimensional vector space 
$\sum_{v,w\in V}v\otimes w+w\otimes v$
of all symmetric tensor products in $V\otimes V$.
In the sequel we speak simply of perfect sets if the ambient 
vector space is obvious.

A choice of a basis $x_1,\dots,x_d$ of $V$ identifies
$V$ with the vector space $\{a_1x_1+\dots+a_dx_d\ \vert\ 
a_1,\dots,a_d\in\mathbb R\}$ of all homogeneous 
$1$-forms in $\mathbb R[x_1,\dots,x_d]$. 
Perfection of $\mathcal S$ is equivalent to the fact that the
set
$$\left\lbrace\left(\sum_{i=1}^d a_i x_i\right)^2\right\rbrace_{\sum_{i=1}^d 
a_i x_i
\in \mathcal S}$$ of all quadratic forms associated to elements in $\mathcal S$
spans the full ${d+1\choose 2}$-dimensional 
vector space of all quadratic forms (homogeneous polynomials of degree $2$).
Equivalently, $\mathcal S$ is perfect (in $V$) if and only if 
the set of symmetric matrices
$\left\{(a_ia_j)_{1\leq i,j\leq d}\right\}_{\sum a_ix_i\in \mathcal S}$
spans the vector space of all symmetric square-matrices of size $d$.
The matrix $(a_ia_j)_{1\leq i,j\leq d}$ is, up to a scalar multiple, 
the orthogonal projection of $V$ onto $\sum a_ix_i$ with respect to 
the scalar product with orthonormal basis $x_1,\dots,x_n$.

Perfect sets of vector spaces over real fields 
determine scalar products uniquely in the following way: 
A scalar product $\langle \ ,\ \rangle:
V\times V\longrightarrow \mathbb R$ on $V\times V$ is uniquely defined
by the set $\{\langle v,v\rangle\}_{v\in\mathcal S}$ of norms of elements in 
$\mathcal S$ if and only if $\mathcal S$ is perfect.

\begin{rem}\label{remperfgen}
Perfection can be generalized as follows: A subset $\mathcal S$ of
a vector space $V$ over a field of characteristic $0$ or larger than $k$
is \emph{$k$-perfect in $V$} 
if the elements of the set 
$\{v\otimes v\otimes \cdots \otimes v\in V^{\otimes^k}\ \vert\ v\in 
\mathcal S\}$ span the ${n+k-1\choose k}$-dimensional
subspace spanned by all symmetric $k$-fold
tensor powers $\{\sum_{\sigma\in S_k}v_{\sigma(1)}\otimes v_{\sigma(2)}\otimes
\cdots\otimes v_{\sigma(k)}\ \vert\ v_1,\dots,v_k\in V\}$ of $V$. A 
$k$-perfect subset of $V$ is $k'$-perfect for
$k'\leq k$. Given a subset $\mathcal S$ of a $d$-dimensional vector
space $V$ over a field $\mathbb K$ of characteristic $0$, we denote by 
$A_k$ the vector space spanned by the set 
$\{v\otimes v\otimes \cdots \otimes v\in V^{\otimes^k}\ \vert\ v\in 
\mathcal S\}$ and by $\alpha_k=\dim(A_k)$ the dimension of $A_k$.
We use the convention $A_0=\mathbb K$ and $\alpha_0=1$.
The generating series $\sum_{k=0}^\infty \alpha_kt^k\in\mathbb N[t]$ 
is always a rational function of the form 
$\frac{P_{\mathcal S}(t)}{1-t}$ for $P_{\mathcal S}\in\mathbb N[t]$
a polynomial with non-negative integral coefficients.
It would be interesting to understand all possible polynomials
arising in this way. In dimension $d=2$ we have
$P_{\mathcal S}=1+t+t^2+\dots+t^{a-1}=\frac{1-t^a}{1-t}$ 
where $a$ is the number of distinct lines 
$\left\{\mathbb Kv\right\}_{v\in \mathcal S}$ 
defined by all elements of $\mathcal S$. 
\end{rem}

We will make repeated use of the following trivial observation
which is assertion 1 of Proposition 3.5.3 in \cite{eM}:

\begin{prop}\label{propfondperf}
Let $\mathcal S$ be a set of non-zero elements in a $d$-dimensional
vector space $V$. Suppose that $V$ contains 
a hyperplane $\mathcal H$ such that 
$\mathcal S\cap \mathcal H$ is perfect in $\mathcal H$
and suppose that the elements 
$\mathcal S\setminus(\mathcal S\cap \mathcal H)$ of $\mathcal S$
in the complement $V\setminus \mathcal H$ of 
$\mathcal H$ generate $V$. Then $\mathcal S$ is perfect.
\end{prop}

We apply Proposition \ref{propfondperf} always in the case where
$V$ is a Euclidean vector space. The hyperplane $\mathcal H$ can 
then be described as the orthogonal subspace $\mathcal H=v^\perp$ 
of a non-zero element $v$ in $V$.

\noindent{\bf Proof of Proposition \ref{propfondperf}} 
We extend a basis $b_1,\dots,b_{d-1}$ of $\mathcal H$ to a basis
$b_1,\dots,b_d$ of $V$. The vector space spanned by 
$\{v\otimes v\}_{v\in\mathcal S\cap \mathcal H}$ contains 
the vector space of all symmetric tensor products in $\mathcal H\otimes
\mathcal H$ by perfection of $\mathcal S\cap \mathcal H$.
The fact that $\mathcal S\setminus(\mathcal S\cap \mathcal H)$
generates $V$ implies that 
the vector space spanned by $\{v\otimes v\}_{v\in\mathcal S\setminus(\mathcal S
\cap \mathcal H)}$ contains all symmetric tensors 
$b_i\otimes b_d+b_d\otimes b_i$ for $i=1,\dots,d$.\hfill$\Box$

\subsection{Perfect lattices}

A \emph{Euclidean lattice} (or \emph{lattice} in the sequel) is
a discrete subgroup of a finite-dimensional Euclidean vector space
$\mathbb E$.
A lattice $\Lambda$ of rank $d=\mathrm{dim}(\Lambda\otimes_\mathbb Z\mathbb R)$ 
is isomorphic to $\mathbb Z^d$ as a group.
$\Lambda$ is \emph{integral} if the scalar 
product $\mathbb E\times \mathbb E\ni(u,v)\longmapsto \langle u,v\rangle
\in\mathbb R$ has an integral restriction 
$\langle \ ,\ \rangle:\Lambda\times \Lambda\longrightarrow \mathbb Z$.
The \emph{norm} of a lattice-element $\lambda$ is in the sequel 
always the squared Euclidean length $\langle \lambda,\lambda\rangle$ of 
$\lambda$. An integral lattice is \emph{even} if all its elements
have even norm.
We denote by $\Lambda_{\min}$ the set of shortest non-zero elements,
called \emph{minimal elements}, in 
$\Lambda$ and by $\min(\Lambda)$ the \emph{minimal norm} $\langle
v,v\rangle$ of a minimal element $v$ in $\Lambda_{\min}$.
The \emph{determinant} $\det(\Lambda)$ of a lattice is 
the squared volume of a fundamental domain for the action (by translations) 
of $\Lambda$ on $(\Lambda\otimes_{\mathbb Z}
\mathbb R)$. The determinant $\det(\Lambda)$ 
is given by $\det(G)$ with $(G)_{i,j}=\langle b_i,b_j)$
a Gram matrix defined by scalar products between basis elements 
$b_1,\dots,b_d$ of the $d$-dimensional lattice $\Lambda=\oplus_{i=1}^d
\mathbb Zb_i$.
The \emph{density} $$\frac{\sqrt{\min(\Lambda)}^d}{2^d\sqrt{\det(\Lambda)}}
\frac{\pi^{d/2}}{(d/2)!}$$ 
of a $d$-dimensional lattice $\Lambda$ is 
the density of the associated sphere-packing obtained
by packing the space $\Lambda\otimes_{\mathbb Z}\mathbb R$ with
spheres of equal radius $\sqrt{\min(\Lambda)}/2$ (and delimiting balls
of volume $\left(\frac{\sqrt{\min(\Lambda)}}{2}\right)^d
\frac{\pi^{d/2}}{(d/2)!}$)
centered at all lattice points. 
\emph{Extreme lattices} are lattices whose density is locally maximal 
(with respect to the obvious natural topology
on the space of lattices of given dimension).
Extreme lattices are perfect and eutactic (a positivity condition), 
cf. Theorem 3.4.6 in \cite{eM}.
Perfection and eutaxy are however independent in the sense that
one property does not necessarily imply the other.
Thus there exist perfect lattices which are not extreme.
All perfect lattices can be realized, up to similarity, as
integral lattices (cf. Proposition 3.2.11 of \cite{eM}) 
and there are only finitely many of them
(up to similarity and isometry) in any given dimension, cf. 
Theorem 3.5.4 in \cite{eM}.  
The following definition provides a measure for perfection:
Given a lattice $\Lambda$ of rank $d$, we denote by 
$\mathrm{pd}(\Lambda)$ its \emph{perfection-default} (called 
co-rank in the monograph \cite{eM} devoted to perfect lattices)
defined as ${d+1\choose 2}-\dim(\mathcal A)$ with $\mathcal A
=\sum_{v\in\Lambda_{\min}}\mathbb R\ v\otimes v$ denoting
the vector space spanned by $\{v\otimes v\}_{v\in\Lambda_{\min}}$.
A lattice is perfect if and only if its perfection-default is zero. 

The aim of this paper is the construction
of a few integral lattices with minimum $4$ (we describe also a 
family with minimum $3$). All considered lattices are sublattices
of $\mathbb Z^n$ and are thus kernels of morphisms
$\varphi:\mathbb Z^n\longrightarrow A$ onto a suitable abelian group $A$.
The specific form of $\varphi$ is of crucial importance
since it allows the deduction of perfection from combinatorial properties.
Our construction is very flexible and gives rise to many 
inequivalent perfect lattices. In particular, we show in Theorem
\ref{thmgrowth} that the number
of inequivalent perfect integral lattices of minimum $4$ and dimension $d$
has no polynomial upper bound as a function of $d$.

The sequel of this paper is organized as follows:
Section \ref{sectperflatt} describes the main construction and its
generalization, obtained by considering suitable $d$-dimensional 
sublattices of the $(d+1)$-dimensional root lattice of type $A$. 
The rest of the paper is essentially a variation on this theme.
Section \ref{sectoddconstr} avoids the use of the root lattice
of type $A$ by considering the orthogonal of an integral vector
having only odd coefficients.
Section \ref{sectevenconstr} replaces the root lattice of type $A$
by the root lattice of type $D$.
Section \ref{sectabelconstr} considers sublattices of finite index 
in root lattices of type $A$.
Section \ref{sectevanabelconstr} considers sublattices of finite index 
in root lattices of type $D$.
Section \ref{sectmin3} discusses briefly a family of perfect lattices 
having minimum $3$ related to projective spaces over the field
$\mathbb F_2$ of $2$ elements.
The rest of the paper deals with other variations
based on finite abelian groups and generalizations.

\section{A sequence of perfect lattices}\label{sectperflatt}

We denote by $L_d$ the even integral lattice of rank $d$ 
defined by all vectors
of $\mathbb Z^{d+2}$ orthogonal to both elements $(1,1,\dots,1)$ and 
$(1,2,\dots,d+2)$ of $\mathbb Z^{d+2}$.

\begin{thm}\label{thmperfectionLdfamily} 
The lattice $L_d$ has determinant $\frac{1}{12}(d+1)(d+2)^2(d+3)$
and contains no roots (vectors of norm $2$). It has 
$\frac{1}{24}d(d+2)(2d-1)$ pairs of opposite vectors of 
(squared Euclidean) norm $4$ if $d$ is even 
and $\frac{1}{24}(d-1)(d+1)(2d+3)$ pairs of opposite vectors of norm $4$
if $d$ is odd. The lattice $L_d$ is perfect for $d\geq 7$.
\end{thm}

\begin{rem} The lattice $L_6$ has $22$ pairs of minimal 
elements. The set $\left\lbrace v\otimes v\right\rbrace_{v\in\min(L_6)}$
spans a vector 
space of dimension $20$. The lattice $L_6$ has thus perfection-default
${7\choose 2}-20=1$ and is not perfect.

The seven rows of the matrix
$$A=\left(\begin{array}{rrrrrrrrr}
1&-1&0&0&0&-1&1&0&0\\
0&0&0&1&-1&-1&1&0&0\\
1&0&-1&0&-1&0&1&0&0\\
0&1&-1&0&0&-1&1&0&0\\
0&-1&1&0&0&0&1&-1&0\\
1&-1&-1&1&0&0&0&0&0\\
1&-1&0&0&0&0&0&-1&1
\end{array}\right)$$
span the perfect lattice $L_7$. The associated Gram matrix $AA^t$ 
with determinant $\frac{2^3\cdot 3^4\cdot 2\cdot 5}{2^2\cdot 3}=
2^2\cdot 3^3\cdot 5$
is the reduced Gram matrix
$$P_7^7=\left(\begin{array}{rrrrrrr}
4&2&2&1&2&2&2\\
2&4&2&2&1&1&0\\
2&2&4&2&0&2&1\\
1&2&2&4&-1&0&-1\\
2&1&0&-1&4&0&2\\
2&1&2&0&0&4&2\\
2&0&1&-1&2&2&4
\end{array}\right)$$
at page 382 of \cite{M}. Gram-matrices for all perfect lattices
up to dimension $7$ are only given in the original French version \cite{M}.
They are unfortunately missing in the English translation \cite{eM}.
\end{rem}

\noindent{\bf Proof of Theorem \ref{thmperfectionLdfamily}}
The determinant of $L_d$ is equal to the determinant of the 
$2$-dimensional lattice $\mathbb Z^{d+2}\cap
(L_d\otimes_{\mathbb Z}\mathbb R)^\perp$ of $\mathbb Z^{d+2}$ 
which is orthogonal to $L_d$.
Thus the determinant of $L_d$ is given by
$$\det\left(\begin{array}{cc}\langle u,u\rangle&\langle u,v\rangle\\
 \langle v,u\rangle&\langle v,v\rangle\end{array}\right)$$
with $u=(1,1,\dots,1),v=(1,2,\dots,d+2)\in\mathbb Z^{d+2}$,
and checking the formula is straightforward.

The lattice $L_d$ is obviously integral and even. 
Vectors of norm $2$ in $\mathbb Z^{d+2}$ are of the form 
$\pm e_i\pm e_j$ (with $e_1,\dots,e_{d+2}$ denoting the natural orthonormal
basis of $\mathbb Z^{d+2}$) and are never orthogonal to 
both elements $(1,\dots,1)$ and $(1,\dots,d+2)$ of $\mathbb Z^{d+2}$.

Vectors of norm $4$ in $L_d$ are of the form 
$$e_i-e_{i+\alpha}-e_{i+\alpha+\beta}+e_{i+2\alpha+\beta}$$
with $i\in\{1,\dots,d-1\}$ and $\alpha,\beta$ two natural numbers
greater than $0$ such that 
$i+2\alpha+\beta\leq d+2$. The lattice $L_d$ contains thus 
$$\sum_{i=1}^{d-1}\sum_{\alpha=1}^{\lfloor (d+1-i)/2\rfloor}
d+2-i-2\alpha$$
pairs of minimal vectors. This formula, restricted to even, respectively 
odd, natural integers, defines a polynomial function
of degree $3$. Explicit expressions can be found
by interpolation of four values.

We prove perfection of $L_d$ by induction on $d$. Perfection of the lattice
$L_7$ considered above establishes the result for $d=7$.

The identity 
$(1,2,3,\dots,d+3)-(1,1,1,\dots,1)=(0,1,2,\dots,d+2)$ 
shows that the sublattice of $L_{d+1}$ orthogonal to $(1,0,0,\dots,0)$
is the lattice $L_d$ which is perfect by assumption.
By Proposition \ref{propfondperf} it is enough to show that
the vector space spanned by the set of minimal vectors with
first coordinate non-zero (i.e. with first coordinate $\pm 1$) 
has dimension $d+1$.
We set $u=e_1-e_3-e_4+e_6$ and $v_i=e_1-e_2-e_{i-1}+e_i$ for $i=4,\dots,d+3$.
Consideration of the last index $i$ with non-zero coefficient
of the vector $v_i$ shows linear independency
of the $d$ vectors $v_4,\dots,v_{d+3}$.
Computation of
$$v_5+v_6-u=e_1-2e_2+e_3$$
shows linear independency of $u$ from $v_4,\dots,v_{d+3}$
and ends the proof.\hfill$\Box$


\subsection{A generalization}

To a strictly increasing sequence $1\leq a_1<a_2<\dots<a_k$
of $k$ natural integers $a_1,\dots,a_k$ and an integer $n$ 
we associate the set $\mathcal I_n(a_1,\dots,a_k)$ 
defined by the smallest $n$ elements of
$\{1,2,\dots\}\setminus\{a_1,\dots,a_k\}$.
We consider now the
sequence of lattices $L_d(a_1,\dots,a_k)$ 
consisting of all elements of 
$\mathbb Z^{d+2}$ which are orthogonal to $(1,\dots,1)\in\mathbb Z^{d+2}$ 
and to the vector of $\mathbb Z^{d+2}$ with increasing 
coefficients given by the elements of $\mathcal I_{d+2}(a_1,\dots,a_k)$.
Equivalently, $L_d(a_1,\dots,a_k)$ can be defined for $d$ large enough 
as the
sublattice of $L_{d+k}$ defined by all vectors with zero coefficients
for indices in $\{a_1,a_2,\dots,a_k\}$.

The lattices $L_d(a_1,\dots,a_k)$ and $L_d(d+2+k-a_k,d+2+k-a_{k-1},\dots,
d+2+k-a_1)$ are obviously isomorphic for $d>a_k-k-2$.

Theorem \ref{thmperfectionLdfamily} has the following generalization:

\begin{thm}\label{thmperfectfamgen} 
The lattice $L_d(a_1,\dots,a_k)$ is perfect for $d\geq \max(7,2(k+1)^3-1)$.
\end{thm}

Theorem \ref{thmperfectfamgen} is an easy consequence of the following
two results:

\begin{prop}\label{propauxil} The lattice $L_d(a_1,\dots,a_k)$ with $d\geq 7$ 
is perfect if $a_1\geq 2(k+1)^2+2$.
\end{prop}

\begin{prop}\label{propauxil1} The lattice $L_d(a_1,\dots,a_k)$ 
is perfect if the subset $\mathcal I_d(a_1,\dots,a_k)$ 
defined by the $d+2$ smallest elements of 
$\{1,2,\dots\}\setminus\{a_1,\dots,a_k\}$ contains 
$\max(2(k+1)^2+1,9)$ consecutive elements.
\end{prop}

\noindent{\bf Proof of Theorem \ref{thmperfectfamgen}}
The result holds for $k=0$ by Theorem \ref{thmperfectionLdfamily}.
Removing a non-empty set of $k$ integers from $\{1,2,\dots\}$ leaves (at most)
$k+1$ subsets of consecutive integers. A partition of 
$2(k+1)^3-1+2=2(k+1)^3+1$ elements into (at most) $k+1$ subsets
of consecutive integers contains thus a subset
having at 
least $2(k+1)^2+1$ consecutive elements and the result follows from
Proposition \ref{propauxil1}.\hfill$\Box$

\noindent{\bf Proof of Proposition \ref{propauxil}} Theorem 
\ref{thmperfectionLdfamily} shows that the result holds
for $d\in\{7,\dots, 2(k+1)^2-1\}$. As in the proof of Theorem 
\ref{thmperfectionLdfamily} we use induction on $d$ 
establishing the induction-step through Proposition 
\ref{propfondperf}. The sublattice of 
$L_{d+1}(a_1,\dots,a_k)$ 
consisting of all elements with last coefficient zero is the  
lattice $L_d(a_1,\dots,a_k)$ which is perfect by assumption.
We have thus only to prove that the set of all minimal
vectors in $L_{d+1}(a_1,\dots,a_k)$ 
with last coefficient non-zero spans a vector space of dimension $d+1$.
For simplicity, we work with $\mathcal I_{d+3}(a_1,\dots,a_k)$
as the set of indices for the $d+3$ coefficients of elements in
$L_{d+1}(a_1,\dots,a_k)$.
We denote by $\omega\leq d+3+k$, respectively $\psi$ the largest,
respectively second-largest, element of $\mathcal I_{d+3}
(a_1,\dots,a_k)$. For every index 
$i\in\mathcal I_{d+3}(a_1,\dots,a_k)$ with $i<\psi$ 
we construct a linear combination $u(i)=\sum \alpha_vv$
of minimal elements $v=(v_1,\dots, v_{\omega})$ ending with 
last non-zero coefficient $v_\omega=1$ 
such that $i$ is
the index of the first non-zero coefficient in $u(i)$.
If $i<\omega-2k-3$, there exists an integer 
$j=j(i)$ in $\{1,\dots,k+1\}$ such that both integers $i+j$ and $\omega-j$ 
belong to $\mathcal I_{d+3}(a_1,\dots,a_k)$ and we can take 
$u(i)=e_i-e_{i+j}-e_{\omega-j}+e_\omega$. 
For $i\in \mathcal I_{d+3}(a_1,\dots,a_k)$ such 
that $i\in\{\omega-2k-2,\dots,\psi-1\}$
we set $\alpha=\omega-i$ and $\beta=\omega-\psi$. We have 
$1\leq \beta\leq k+1$ and $\beta\leq \alpha\leq 2(k+1)$.
Since $a_1\geq 2(k+1)^2+2$, all integers $1,2,\dots,1+\alpha\beta\leq
2(k+1)^2+1$ are in $\mathcal I$ and we can consider
\begin{eqnarray*}
u(i)&=&(e_1-e_{1+\alpha}-e_i+e_\omega)+(e_{1+\alpha}-e_{1+2\alpha}-e_i+
e_\omega)\\
&&+\dots+(e_{1+(\beta-1)\alpha}-e_{1+\beta\alpha}-e_i+e_\omega)\\
&&-(e_1-e_{1+\beta}-e_\psi+e_\omega)-(e_{1+\beta}-e_{1+2\beta}-e_\psi+
e_\omega)\\
&&-\dots-(e_{1+(\alpha-1)\beta}-e_{1+\alpha\beta}-e_\psi+e_\omega)\\
&=&-\beta e_i+\alpha e_\psi+(\beta-\alpha)e_\omega
\end{eqnarray*}
which ends the proof.\hfill$\Box$

\noindent{\bf Proof of Proposition \ref{propauxil1}}
Theorem \ref{thmperfectionLdfamily} shows the result for $k=0$.
We assume henceforth $k\geq 1$ and $\max(2(k+1)^2+1,9)=
2(k+1)^2+1$.

We denote by $\alpha$
the smallest integer such that $\{\alpha,\alpha+1,\dots,
\alpha+2(k+1)^2\}$ is contained in $\{1,2,\dots\}\setminus\{a_1,\dots,a_k\}$.
For $\alpha=1$ the proof follows from Proposition \ref{propauxil}.
Otherwise we establish the result by induction on $k$ and $d$ using 
Proposition \ref{propfondperf}. 
We assume $\alpha>1$ and we consider $L_{d+1}(a_1,\dots,a_k)$.
Since $L_{d+1}(a_1,\dots,a_k)$ with $a_1=1$ 
is isomorphic to $L_{d+1}(a_2-1,a_3-1,\dots,a_k-1)$ (which 
is perfect by induction on $k$, the case $k=0$ being covered by 
Theorem \ref{thmperfectionLdfamily}) we can assume $a_1>1$ .
Since $\alpha>1$, the sublattice $L_d(a_1-1,a_2-1,\dots,a_k-1)$ 
consisting of all vectors of $L_{d+1}(a_1,\dots,a_k)$
with first coefficient zero is perfect by induction on $d$.
By Proposition \ref{propfondperf} it is thus enough to show that the set of
minimal elements of $L_{d+1}(a_1,\dots,a_k)$ with first coefficient non-zero
spans a $(d+1)$-dimensional vector space. 
The proof is analogous to the proof of Proposition \ref{propauxil}
except that we work with small coordinates instead of large ones and
that we use the $2(k+1)^2+1$ indices $\alpha,\dots,\alpha+2(k+1)^2$
instead of the set $\{1,\dots,2(k+1)^2+1\}$.\hfill$\Box$

\subsection{Examples for Theorem \ref{thmperfectfamgen}}
Below we list a few lattices of dimension $7$ or $8$ illustrating
Theorem \ref{thmperfectfamgen} for 
$k=1,2$. For $k=1$ we indicate the relevant integer $a_1$ missing in
$1,2,\dots,d+3$ together with the determinant, the perfection-default
$\hbox{pd}(L_d(a_1))={d+1\choose 2}-\dim\left(\sum_{v\in L_d(a_1)}
\mathbb R\ v\otimes v\right)$ and the number $\mathrm{mp}$ of pairs of 
minimal vectors.
The lattice $L_d(1)$ is obviously isomorphic to the lattice 
$L_d=L_d(\emptyset)$ studied previously.
Since $L_d(i)$ and $L_d(d+4-i)$
are isomorphic, it is enough to consider the four cases $i=2,3,4,5$ 
for $d=7$.
$$\begin{array}{|l|c|r|r|}
\hline
\hbox{lattice}&\det&\hbox{pd}&\hbox{mp}\\
\hline
\hline
L_7(2)&2^2\cdot 5\cdot 31&1&31\\
\hline
L_7(3)&2^3\cdot 5\cdot 17&1&29\\
\hline
L_7(4)&2^4\cdot 3^2\cdot 5&0&28\\
\hline
L_7(5)&2^2\cdot 5\cdot 37&4&28\\
\hline
\end{array}$$

The perfect lattice $L_7(4)$ of 
determinant $2^4\cdot 3^2\cdot 5$
(and defined as the set of all integral vectors
in $\mathbb Z^9$ which are orthogonal to 
$(1,1,1,1,1,1,1,1,1)$ and $(1,2,3,5,6,7,8,9,10)$)
has a basis given by the rows of 
$$A=\left(\begin{array}{ccccccccc}
0&1&-1&0&0&0&0&-1&1\\
0&0&-1&1&0&1&0&-1&0\\
0&0&0&1&-1&0&0&-1&1\\
0&0&0&0&0&1&-1&-1&1\\
1&-1&0&0&-1&1&0&0&0\\
1&0&-1&0&0&0&-1&0&1\\
1&-1&0&0&0&0&0&-1&1
\end{array}\right)$$
with Gram matrix $AA^t$ given by the matrix
$$P_7^{31}=\left(\begin{array}{ccccccc}
4&2&2&2&-1&2&1\\
2&4&2&2&1&1&1\\
2&2&4&2&1&1&2\\
2&2&2&4&1&2&2\\
-1&1&1&1&4&1&2\\
2&1&1&2&1&4&2\\
1&1&2&2&2&2&4
\end{array}\right)$$
at page 383 of \cite{M}. 

The list of lattices of the 
form $L_8(a_i)$ (with $a_i\in\{2,\dots,6\}$
in order to avoid duplicates) with a few properties is given by
$$\begin{array}{|l|c|r|r|}
\hline
\hbox{lattice}&\det&\hbox{pd}&\hbox{mp}\\
\hline
\hline
L_8(2)&2^2\cdot 3\cdot 7\cdot 11&0&46\\
\hline
L_8(3)&7\cdot 11\cdot 13&0&44\\
\hline
L_8(4)&2^5\cdot 3\cdot 11&0&42\\
\hline
L_8(5)&3^2\cdot 11^2&0&42\\
\hline
L_8(6)&2^2\cdot 5^2\cdot 11&3&42\\
\hline
\end{array}$$

The following table lists all six perfect 
non-isomorphic lattices of the form $L_8(a_1,a_2)$ obtained by
removing two elements $a_1,a_2$ from 
$\{2,\dots,11\}$ (we exclude $a_1=1$ in order to avoid 
perfect lattices of the form $L_8(a_1)$):
$$\begin{array}{|l|c|r|r|}
\hline
\hbox{lattice}&\det&\hbox{pd}&\hbox{mp}\\
\hline
\hline
L_8(2,3)&3\cdot 347&0&43\\
\hline
L_8(2,5)&7\cdot 167&0&40\\
\hline
L_8(2,6)&2^4\cdot 3\cdot 5^2&0&39\\
\hline
L_8(2,9)&3^3\cdot 43&0&40\\
\hline
L_8(2,10)&2^4\cdot 3\cdot 23&0&41\\
\hline
L_8(3,5)&2^2\cdot 3^2\cdot 5\cdot 7&0&37\\
\hline
\end{array}$$

\subsection{Bounds for perfection}

Every finite sequence $a_1<a_2<\dots<a_k$ determines a maximal subset
$\mathcal P(a_1,\dots,a_k)$ of $\mathbb N$ such that 
$L_d(a_1,\dots,a_k)$ has minimum $4$ and is perfect
for $d\in\mathcal P(a_1,\dots,a_k)$.
We denote by $D(a_1,\dots,a_k)$ the successor of the largest missing integer
in $\mathcal P(a_1,\dots,a_k)$. We have $D(a_1,\dots,a_k)\leq 
\max(7,2(k+1)^3-1)$ by Theorem  \ref{thmperfectfamgen}.
Since $L_d(a_1,\dots,a_k)$ is perfect for every 
$d\geq \max(7,2(k+1)^3-1)\geq D(a_1,\dots,a_k)$  
there exists a 
smallest integer $d_k=\max_{a_1,\dots,a_k}D(a_1,\dots,a_k)$ 
(bounded above by $\max(7,2(k+1)^3-1)$) 
such that $L_d(a_1,\dots,a_k)$ is perfect for every $d\geq d_k$ and
for every $\{a_1,\dots,a_k\}$ in $\{1,2,\dots\}$. We have $d_0=7$
by Theorem \ref{thmperfectionLdfamily}.

For $k=1$ we get the numbers
$$\begin{array}{|c||c|c|c|c|c|c|c|c|c|c||}
\hline
a_1&1&2&3&4&5&6&7&8&9&\geq 10\\
\hline
D(a_1)&7&8&8&7&8&9&7&8&8&7\\
\hline\end{array}$$
showing $d_1=9$.

\begin{rem} Analogues of the above numbers and bounds exist
of course also for most subsequent constructions.
\end{rem} 

\subsection{Automorphisms and growth}

The aim of this Section is to sketch a proof of the following result:

\begin{thm} \label{thmgrowth}
The number of non-isomorphic perfect integral lattices of
dimension $d$ and minimum $4$ grows faster than any polynomial in 
$d$.
\end{thm}

The proof of Theorem \ref{thmgrowth} is based on Theorem 
\ref{thmperfectfamgen} and gives an explicit lower bound on 
the number of perfect integral lattices of dimension $d$ and
minimum $4$. This lower bound is unlikely to be sharp: The construction
underlying Theorem 
\ref{thmperfectfamgen} yields probably only a small fraction of
all non-isomorphic integral perfect lattices with minimum $4$. 
Moreover, the bounds in Theorem 
\ref{thmperfectfamgen} are certainly far from optimal.

Two minimal vectors $v,w\in(L_d)_{\min}$ are \emph{neighbours}
if $\langle v,w\rangle=2$.

We call the real number 
$\gamma(v)=\frac{2(i+\alpha)+\beta}{2(d+1)}\in (0,1)$ the \emph{(normalized)
center}
and $\delta(v)=\frac{2\alpha+\beta}{d+1}\in (0,1]$ the \emph{(normalized)
diameter} of a minimal vector
$$v=\pm(e_i-e_{i+\alpha}-e_{i+\alpha+\beta}+e_{i+2\alpha+\beta})$$
in $L_d$. We have $\delta(v)\leq 2\min(\gamma(v),1-\gamma(v))$.

\begin{lem}\label{lemformnberneighbrs} 
The number of neighbours in $L_d$ of a minimal vector $v$
in $L_d$ is given by
\begin{eqnarray}\label{formvois}
&2d\left(\min(\gamma(v),1-\gamma(v))+2-\delta(v)\right)+O(1)\ .
\end{eqnarray}
\end{lem}

Formula (\ref{formvois}) is bounded above by $5d$ 
with asymptotic equality for 
$\gamma(v)=\frac{1}{2}$ and $\delta(v)=0$ and bounded below by 
$3d$ with asymptotic equality
for $\gamma(v)=\frac{1}{2}$ and $\delta(v)=1$.

\noindent{\bf Proof of Lemma \ref{lemformnberneighbrs}}
Neighbours of a minimal vector $v\in \min(L_d)$ are partitioned into 
$6={4\choose 2}$ families
according to their two common non-zero coefficients.
We denote these families by 
$$\mathcal F_{**00},\mathcal F_{*0*0},\mathcal F_{*00*},\mathcal F_{0**0},
\mathcal F_{0*0*},\mathcal F_{00**}$$ 
where $*$ stands for a common
non-zero coefficient with respect to the obvious linear order $i<i+\alpha<
i+\alpha+\beta<i+2\alpha+\beta$ on the indices of the four coefficients of
$v=e_i-e_{i+\alpha}-e_{i+\alpha+\beta}+e_{i+2\alpha+\beta}$.
Neighbours in $\mathcal F_{*00*}$ or in $\mathcal F_{0**0}$
share the center with $v$ and 
there are roughly $d\min(\gamma(v),1-\gamma(v))$
possibilities for the remaining smallest nonzero coordinate in each of these
families.
 
The number of neighbours in $\mathcal F_{**00}$ (or in $\mathcal F_{00**}$)
is roughly given by $d-\alpha$ and the number of neighbours in
$\mathcal F_{*0*0}$ (or in $\mathcal F_{0*0*}$) is roughly given by 
$d-\alpha-\beta$. 

All errors are bounded by an absolute constant (which is small).
Summing over all families and using the definitions of the normalized
center and diameter ends the proof.\hfill$\Box$

\noindent{\bf Sketch of proof for Theorem \ref{thmgrowth}}
We partition neighbours of a minimal vector $v$ into six families
as in the proof of Lemma \ref{lemformnberneighbrs}.
An element $u$ of such a family $\mathcal F
\in\{\mathcal F_{**00},\mathcal F_{*0*0},\mathcal F_{*00*},\mathcal F_{0**0},
\mathcal F_{0*0*},\mathcal F_{00**}\}$ is adjacent to all other 
elements of $\mathcal F$
except perhaps for two elements $u',u''$ for which we have 
$\langle u,u'\rangle=\langle u,u''\rangle=-1$. In $\mathcal F_{*00*}$
and $\mathcal F_{0**0}$ there are no exceptions: two distinct elements of 
$\mathcal F_{*00*}$ or of $\mathcal F_{0**0}$ are always adjacent.
Families associated to complementary pairs of indices, like
$\mathcal F_{*00*},\mathcal F_{0**0}$ or $\mathcal F_{**00},\mathcal F_{00**}$ 
or the remaining two sets $\mathcal F_{0*0*},\mathcal F_{*0*0}$ are
called \emph{complementary}. Complementary sets are related by
a natural involution $\iota$ defined by $\iota(u)=u'$
for $u\in\mathcal F$ and $u'\in\overline{\mathcal F}$
such that $\langle u,u'\rangle =-2$.
Complementary sets have the same number
of elements, given by $d\min(\gamma(v),1-\gamma(v))+O(1)$ for 
$\mathcal F_{*00*}$ or for $\mathcal F_{0**0}$.
The sets $\mathcal F_{**00},\mathcal F_{00**}$ associated to 
$v=e_i-e_{i+\alpha}-e_{i+\alpha+\beta}+e_{i+2\alpha+\beta}$ have 
$d-\alpha+O(1)$ elements and the sets
$\mathcal F_{0*0*},\mathcal F_{*0*0}$ have $d-\alpha-\beta+O(1)$ elements.
An element $u\in \mathcal F$ is orthogonal to every element $v$ of the 
complementary pair $\overline{\mathcal F}$, except for $\iota(u)$
and perhaps at most
two other elements in $\overline{\mathcal F}$. The behaviour is simpler in
the complementary pair $\mathcal F_{*00*},\mathcal F_{0**0}$:
every element $u$ of $\mathcal F_{*00*}$ is orthogonal to every 
element of $\mathcal F_{0**0}\setminus\{\iota(u)\}$. 

These properties allow to reconstruct (at least approximately,
however exact coordinates can be found) 
the coordinates, up to the obvious symmetry $e_i\longmapsto e_{d+3-i}$,  
of a minimal element $v$
from the knowledge of all scalar products between elements of $(L_d)_{\min}$.
This shows that the lattice $L_d$ has at most four automorphisms:
$\pm 1$, perhaps followed by reversal of all coordinates.
The same holds for the lattices $L_d(a_1,\dots,a_k)$: Knowledge 
of all scalar products between the set $(L_d(a_1,\dots,a_k))_{\min}$
of minimal elements determines, up to signs and
global reversal of all coordinates,
the coefficients (and thus also the missing integers $a_1,\dots,a_k$).
In particular, all lattices $L_d(a_1,\dots,a_k)$ have at most four 
automorphisms. In the generic case, only the two trivial automorphisms
$\pm 1$ occur. Thus the number of different lattices
$L_d(a_1,\dots,a_k)$ of large dimension $d$ is at least equal to 
$\frac{1}{2}{d+2+k\choose k}\geq \frac{d^k}{2\cdot k!}$. 
Since $k$ is arbitrary, the number of perfect integral lattices
of dimension $d$ with minimum $4$ grows faster than 
any polynomial in $d$.\hfill$\Box$

\section{The odd construction}\label{sectoddconstr}

We denote by $O_d$ the $d$-dimensional lattice of all integral
vectors in $\mathbb Z^{d+1}$ which are orthogonal to 
$(1,3,5,\dots,2d+1)$. The lattice $O_d$ is even and contains
no roots. 
Since $(1,3,5,7,\dots)=2(1,2,3,4,\dots)-(1,1,1,1,\dots)$ the lattice 
$O_d$ contains 
$L_{d-1}$ as a sublattice. Pairs of minimal vectors of $O_d$
not contained in the sublattice $L_{d-1}$ are of the form
$e_{2a-1}+e_{2b-1}+e_{2c-1}-e_{2k-1}$ (with indices given by the coefficients
of $(1,3,5,\dots)$) corresponding to 
sums $2k-1=(2a-1)+(2b-1)+(2c-1)$ 
of three distinct odd natural integers $2a-1,2b-1,2c-1\in\{1,3,5,\dots,
2d-3\}$ adding up to $2k-1\in\{9,11,\dots,2d+1\}$.

\begin{thm}\label{thmoddperfect} The lattice $O_d$ has determinant $\frac{1}{3}(d+1)(2d+1)(2d+3)$
and minimum $4$ (for $d\geq 3$). It has 
$$\frac{1}{18}\left(2d^3-3d^2-3d+c_d\right)$$
pairs of minimal vectors of norm 4 where 
$$c_d=\left\lbrace\begin{array}{ll}
0\quad&\hbox{if }d\equiv 0\pmod 3,\\
4&\hbox{if }d\equiv 1\pmod 3,\\
2&\hbox{if }d\equiv 2\pmod 3.\end{array}\right.$$
The lattice $O_d$ is perfect for $d\geq 8$.
\end{thm}


The lattice $O_7$ (with $29$ pairs of minimal vectors) 
has perfection default $1$ and is thus not perfect.

\noindent{\bf Proof of Theorem \ref{thmoddperfect}}
The squared Euclidean norm of $(1,3,5,\dots,2d+1)$
is polynomial of degree $3$ in $d$ and the formula for the 
determinant of $O_d$ (given by $1^2+3^2+5^2+\dots+(2d+1)^2$) 
can thus be checked using $4$ values.

The lattice $O_d$ cannot contain roots and its minimum 
is obviously $4$ (realized e.g. by the vector $(1,-1,-1,1,0,0,\dots)$)
if $d\geq 3$.
The number of pairs of minimal vectors are polynomial functions
of degree $3$ for $d$ in arithmetic progressions of length $6$. 
Computation of small examples gives enumerative formulae.

For proving perfection we use Proposition \ref{propfondperf}
with $\mathcal H=(1,1,1,\dots,1)^\perp\subset\mathbb R^{d+1}$. The inclusion 
of $L_{d-1}=\mathcal H\cap O_d$ in $O_d$ and perfection of $L_{\geq 7}$
(see Theorem \ref{thmperfectionLdfamily})
implies that it is enough to show that minimal 
vectors with non-zero coordinate-sum span a $d$-dimensional 
vector space. As already mentioned, such minimal elements are of the form 
$\pm(e_a+e_b+e_c-e_k)$ with $(2a-1)+(2b-1)+(2c-1)=(2k-1)$ for four
distinct elements $a,b,c,k$ in $\{1,\dots,d+1\}$. The seven minimal
vectors given by the rows of
$$\left(\begin{array}{cccccccc}
1&1&1&0&-1&0&0&0\\
1&1&0&1&0&-1&0&0\\
1&1&0&0&1&0&-1&0\\
1&0&1&1&0&0&-1&0\\
1&1&0&0&0&1&0&-1\\
1&0&1&0&1&0&0&-1\\
0&1&1&1&0&0&0&-1
\end{array}\right)$$
of $O_7$ are linearly independent and span thus the full 
$7$-dimensional vector space orthogonal to
$(1,3,5,7,9,11,13,15)$. The union of these vectors (extended to 
elements of $\mathbb Z^{d+1}$ by appending
zeros) with vectors 
\begin{eqnarray*}
&(1,1,0,0,0,0,1,0,-1,0,0,0,0,0,\dots)\\
&(1,1,0,0,0,0,0,1,0,-1,0,0,0,0,\dots)\\
&(1,1,0,0,0,0,0,0,1,0,-1,0,0,0,\dots)\\
&\quad \vdots\end{eqnarray*}
is a basis of the $d$-dimensional vector space 
$O_d\otimes_{\mathbb Z}\mathbb R=(1,3,\dots,2d+1)^\perp$.
\hfill$\Box$

For a finite increasing sequence $1\leq a_1<a_2<\dots<a_k$
of $k$ odd natural integers, we denote by 
$O_d(a_1,\dots,a_k)\subset\mathbb Z^{d+1}$ 
the $d$-dimensional lattice of all integral vectors
orthogonal to $(1,3,\dots,a_1-2,\widehat{a_1},a_1+2,\dots,
a_{k}-2,\widehat{a_k},a_k+2,\dots,2(d+k)+1)$ 
(elements $\widehat{a_k}$ carrying a magical hat
are removed) with increasing coefficients
given by the $d+1$ smallest elements of $\{1,3,5,7,\dots\}\setminus
\{a_1,\dots,a_k\}$.

The following analogue of Theorem \ref{thmperfectfamgen} holds:

\begin{thm}\label{thmoddperfectfamgen} Given a strictly increasing sequence
$1\leq a_1<a_2<\dots<a_k$ of $k$ odd natural integers,
the $d$-dimensional lattice $O_d(a_1,\dots,a_k)$ 
is perfect for $d\geq \max(10(k+1)^3+5,22(k+1)+2)$.
\end{thm}

\noindent{\bf Proof} 
We apply again Proposition \ref{propfondperf}. Since $O_d(a_1,\dots,a_k)$
contains
the lattice $L_{d-1}((a_1+1)/2,(a_2+1)/2,\dots,(a_k+1)/2))$ defined by all 
elements of $O_d(a_1,\dots,a_k)$ which are orthogonal to 
$(1,1,\dots,1)$ and to 
$\frac{1}{2}\left((1,1,\dots,1)+(1,3,5,\dots,\widehat{a_1},\dots)\right)$, 
we suppose $d$ large enough in order to ensure perfection of
the sublattice $L_{d-1}((a+1)/2,(a_2+1)/2,\dots,(a_k+1)/2))$ (which can be done
using Theorem \ref{thmperfectfamgen}).
We show now that minimal vectors $e_{2a+1}+e_{2b+1}+e_{2c+1}-e_{2l+1}$
(with indices in $\{1,3,5,\dots\}\setminus\{a_1,\dots,a_k\}$) 
of $O_d(a_1,\dots,a_k)$ corresponding to sums 
$2l+1=2a+1+2b+1+2c+1$ with $2l+1,2a+1,2b+1,2c+1\not\in \{a_1,\dots,a_k\}$)
span the $d$-dimensional vector space orthogonal to
$(1,3,\dots,a_1-2,\widehat{a_1},a_1+2,\dots,
a_{k}-2,\widehat{a_k},a_k+2,\dots,2(d+k)+1)$.
This is done in two steps. First we show that for every possible index
$2l+1\geq 12(k+1)+1$ occurring in elements of $O_d(a_1,\dots,a_k)$ 
there exists a minimal vector $e_{2a+1}+e_{2b+1}+e_{2c+1}-e_{2l+1}$
in $O_d(a_1,\dots,a_k)$. The second step 
deals with the remaining vectors involving only small indices.

Every even integer $2m\geq 4(k+1)$ can be written in $\lfloor m/2\rfloor
\geq k+1$
different ways as a sum of two odd natural numbers. For any odd integer
$2l+1\geq 12(k+1)+1$ there exist thus three distinct odd natural numbers 
$2a+1,2b+1,2c+1\not\in\{a_1,\dots,a_k\}$ with $2a+1,2b+1<2c+1$ and 
$2c+1\in\{2l+1-6(k+1),2l+1-4(k+1)\}$ such that $2l+1=2a+1+2b+1+2c+1$:
Indeed, we start by choosing an odd integer 
$2c+1\in\{2l+1-6(k+1),2l+1-4(k+1)\}\setminus\{a_1,\dots,a_k\}$.
The even integer $2l+1-(2c+1)\in\{4(k+1),\dots,6(k+1)\}$
can now be written as a sum of two distinct odd integers $2a+1,2b+1$
not in $\{a_1,\dots,a_k\}$. Since $2l+1\geq 12(k+1)+1$, we have 
$2c+1>\frac{2l+1}{2}$ and thus all three integers $2a+1,2b+1,2c+1
\in\{1,3,5,\dots,\}\setminus\{a_1,\dots,a_k\}$
are distinct. This completes the proof of the first step.

We end now the proof by showing that the span of minimal vectors 
as above (i.e. of the form
$e_{2a+1}+e_{2b+1}+e_{2c+1}-e_{2l+1}$) with $2l+1\geq 12(k+1)+1$ contains the 
vector space $V$ of all vectors which are orthogonal to 
$(1,3,\dots,\widehat{a_1},\dots)$
and which involve only non-zero coefficients
with indices $\leq A=\max(6(k+1)^3+3,14(k+1)+1)$. We have $2a+1+2b+1\leq 2A-2$ for 
two distinct odd integers $a,b\leq A$. There exists thus an 
odd integer $2c+1\in \{A+3,\dots,A+2(k+1)+3\}$
such that $2c+1$ and $2l+1=2a+1+2b+1+2c+1\leq 3A+2(k+1)+1\leq 
\max(20(k+1)^3+10,44(k+1)+4)
\leq 2d+1$ 
are different 
from $a_1,\dots,a_k$. Given four distinct odd integers $2a+1,2b+1,2\alpha+1,
2\beta+1$
in $\{1,3,\dots,A\}\setminus\{a_1,a_2,\dots,a_k\}$ such that 
$2a+1+2b+1=2\alpha+1+2\beta+1$, we can consider the vector
\begin{eqnarray*}
&&(e_{2a+1}+e_{2b+1}+e_{2c+1}-e_{2(a+b+c)+3})\\
&&-(e_{2\alpha+1}+e_{2\beta+1}+e_{2c+1}-e_{2(a+b+c)+3})\\
&=&e_{2a+1}-e_{2\alpha+1}-e_{2\beta+1}+e_{2b+1}.
\end{eqnarray*}
Since the vector space orthogonal to $(1,3,\dots,\hat{a_1},\dots)$
involving no indices exceeding $A$ is at least of dimension 
$$\frac{A-1}{2}-k-1\geq \frac{6(k+1)^3+2}{2}-k-1\geq 2(k+1)^3,$$
minimal
elements of the form $e_{2a+1}-e_{2\alpha+1}-e_{2\beta+1}+e_{2b+1}$
with all indices $\leq A$ span a perfect lattice
by Theorem \ref{thmperfectfamgen}. 
Adding a vector of the form $e_{2a+1}+e_{2b+1}+e_{2c+1}-e_{2l+1}$ with
$2l+1\in\{12(k+1)+1,\dots,14(k+1)+1\}$ which exists by the
discussion of step 1 we get a generating set of $V$.
This completes the proof of the second step and establishes 
Theorem \ref{thmoddperfectfamgen}.\hfill$\Box$

\begin{rem} The bound in Theorem \ref{thmoddperfectfamgen}
(and similar bounds occurring elsewhere) is not optimal
and can be improved by more careful arguments.
\end{rem}

A few data for the lattices 
$O_8(i)$ with $\hbox{pd}$ indicating the perfection default 
and with $\hbox{mp}$ 
indicating the number of pairs of minimal vectors are
$$\begin{array}{|l|c|r|r|}
\hline
\hbox{lattice}&\det&\hbox{pd}&\hbox{mp}\\
\hline
\hline
O_8(1)&3\cdot 443&4&38\\
\hline
O_8(3)&1321&3&37\\
\hline
O_8(5)&3^2\cdot 5\cdot 29&1&37\\
\hline
O_8(7)&3\cdot 7\cdot 61&0&38\\
\hline
O_8(9)&1249&2&38\\
\hline
O_8(11)&3\cdot 13\cdot 31&1&39\\
\hline
O_8(13)&3^3\cdot 43&0&40\\
\hline
O_8(15)&5\cdot 13\cdot 17&1&41\\
\hline
O_8(17)&3\cdot 347&0&43\\
\hline
\end{array}$$
For the lattices $O_9(i)$ the data are
$$\begin{array}{|c|c|r|r|}
\hline
\hbox{lattice}&\det&\hbox{pd}&\hbox{mp}\\
\hline
\hline
O_9(1)&2\cdot 3\cdot 5\cdot 59&2&59\\
\hline
O_9(3)&2\cdot 881&0&56\\
\hline
O_9(5)&2\cdot 3^2\cdot 97&0&56\\
\hline
O_9(7)&2\cdot 3\cdot 7\cdot 41&0&56\\
\hline
O_9(9)&2\cdot 5\cdot 13^2&0&57\\
\hline
O_9(11)&2\cdot 3\cdot 5^2\cdot 11&0&58\\
\hline
O_9(13)&2\cdot 3^2\cdot 89&0&59\\
\hline
O_9(15)&2\cdot 773&0&60\\
\hline
O_9(17)&2\cdot 3\cdot 13\cdot 19&0&62\\
\hline
O_9(19)&2\cdot 3\cdot 5\cdot 47&0&64\\
\hline
\end{array}$$
The lattice $O_{10}(1)$ (with determinant $11^2\cdot 19$ and $81$ pairs
of minimal vectors) is perfect.

\section{The even-sublattice construction}\label{sectevenconstr}

The even-sublattice construction is defined as the $d$-dimensional
lattice $M_d$ consisting of all integral vectors 
$(x_0,\dots,x_d)\in\mathbb Z^{d+1}$
which are orthogonal to $(0,1,2,\dots,d)$ and have
even coordinate-sum $\sum_{i=0}^d x_i\equiv 0\pmod 2$.
Its minimal vectors are $(2,0,0,\dots)$, vectors 
$\pm \left(e_i-e_{i+\alpha}-e_{i+\alpha+\beta}+e_{i+2\alpha+\beta}\right)$ with 
$i\in\{0,\dots,d\}$ and $\alpha,\beta\geq 1$ such that 
$i+2\alpha+\beta\leq d$ together with vectors
$\pm\left(e_h+e_i+e_j-e_k\right)$ where $h,i,j,k$ are four 
distinct integers in $\{0,\dots,d\}$ such that $h+i+j=k$.

\begin{thm}\label{thmevensublconstr}
The lattice $M_d$ has determinant $\frac{2}{3}d(d+1)(2d+1)$
and minimum $4$. It has
$$\frac{1}{36}\left(4d^3-3d^2-6d+c_d\right)$$
pairs of minimal vectors of norm $4$ where $c_d$ depends only on
$d\pmod 6$ and is given by
$$\begin{array}{|c||c|c|c|c|c|c|}
\hline d\pmod 6&0&1&2&3&4&5\\
\hline c_d&36&41&28&45&32&37\\
\hline\end{array}.$$
The lattice $M_d$ is perfect for $d\geq 8$.
\end{thm}

The vector space $\sum_{v\in(M_7)_{\min}}\mathbb R\ v\otimes v$
associated to the $34$ pairs of minimal vectors in $M_7$ is 
of dimension $27$.
The lattice $M_7$ with perfection-default $28-27=1$ is thus not perfect.

\noindent{\bf Proof of Theorem \ref{thmevensublconstr}}
$M_d$ contains $L_{d-1}$ as a sublattice.
Since $d\geq 8$, the lattice $L_{d-1}$ is perfect by Theorem 
\ref{thmperfectionLdfamily}. 
Hence Proposition \ref{propfondperf} shows that it is enough to 
prove that minimal vectors with non-zero coordinate-sum
span the full $d$-dimensional space $M_d\otimes_{\mathbb Z}\mathbb R=
(0,1,2,\dots,d)^\perp\subset \mathbb R^{d+1}$.
This holds for $d\geq 5$ by 
linear independency of the five rows
$$\begin{array}{cccccc}
2&0&0&0&0&0\\
-1&-1&-1&1&0&0\\
-1&-1&0&-1&1&0\\
-1&-1&0&0&-1&1\\
-1&0&-1&-1&0&1\end{array}$$
(with a suitable number of additional zero-coordinates)
together with minimal elements of the form 
$-e_0-e_1-e_{i-1}+e_{i}$ for $i=6,7,\dots,d$.\hfill$\Box$

For a strictly increasing sequence $0\leq a_1<a_2<\dots<a_k$
of $k$ natural integers, we denote by 
$M_d(a_1,\dots,a_k)$ the $d$-dimensional lattice of all integral vectors
with even coordinate-sum which are 
orthogonal to $(0,1,\dots,a_1-1,\widehat{a_1},a_1+1,\dots,
a_{k}-1,\widehat{a_k},a_k+1,\dots,d+k)\in\mathbb Z^{d+1}$.

The following analogue of Theorem \ref{thmevensublconstr} holds:

\begin{thm}\label{thmevensublconstrgen} Given a finite 
strictly increasing sequence
$0\leq a_1<a_2<\dots<a_k$ of $k$ natural integers,
the $d$-dimensional lattice $M_d(a_1,\dots,a_k)$ 
is perfect for $d$ large enough.
\end{thm}
The proof, similar to the proof of Theorem \ref{thmoddperfectfamgen}, is
left to the reader.

A few data with $\hbox{pd}$ indicating the perfection default 
and with $\hbox{mp}$ indicating the number of minimal pairs in lattices 
$M_8(i)$ are
$$\begin{array}{|l|c|r|r|}
\hline
\hbox{lattice}&\det&\hbox{pd}&\hbox{mp}\\
\hline
\hline
M_8(0)&2^2\cdot 3\cdot 5\cdot 19&1&41\\
\hline
M_8(1)&2^4\cdot 71&1&42\\
\hline
M_8(2)&2^2\cdot 281&1&42\\
\hline
M_8(3)&2^4\cdot 3\cdot 23&0&42\\
\hline
M_8(4)&2^2\cdot 269&3&44\\
\hline
M_8(5)&2^4\cdot 5\cdot 13&2&44\\
\hline
M_8(6)&2^2\cdot 3\cdot 83&0&45\\
\hline
M_8(7)&2^4\cdot 59&0&47\\
\hline
M_8(8)&2^2\cdot 13\cdot 17&1&49\\
\hline
\end{array}$$
For the lattices $M_9(i)$ the data are
$$\begin{array}{|c|c|r|r|}
\hline
\hbox{lattice}&\det&\hbox{pd}&\hbox{mp}\\
\hline
\hline
M_9(0)&2^2\cdot 5\cdot 7\cdot 11&0&61\\
\hline
M_9(1)&2^9\cdot 3&0&61\\
\hline
M_9(2)&2^2\cdot 3\cdot 127&0&62\\
\hline
M_9(3)&2^5\cdot 47&0&61\\
\hline
M_9(4)&2^2\cdot 3^2\cdot 41&0&64\\
\hline
M_9(5)&2^5\cdot 3^2\cdot 5&1&64\\
\hline
M_9(6)&2^2\cdot 349&0&65\\
\hline
M_9(7)&2^6\cdot 3\cdot 7&0&66\\
\hline
M_9(8)&2^2\cdot 3\cdot 107&0&69\\
\hline
M_9(9)&2^6\cdot 19&0&70\\
\hline
\end{array}$$

\section{A construction using finite abelian groups}
\label{sectabelconstr}

To a finite abelian group $A$ indexing the coordinates of 
$\mathbb Z^A$ we associate the integral lattice $L(A)$ 
consisting of all elements $v=(v_a)_{a\in A}\in\mathbb Z^{A}$
such that $\sum_{a\in A}v_a=0\in \mathbb Z$ and $\sum_{a\in A}
v_aa=0\in A$ (i.e. vectors $v\in\mathbb Z^A$ of coefficient-sum zero such that 
the element $\sum_{a}v_aa$ of $A$ is the identity $0$ of the finite 
additive group $A$). Equivalently, $L(A)$ is the set of all 
elements in the kernel of the augmentation ideal in the group-algebra
$\mathbb Z[A]$ of $A$ over $\mathbb Z$.
The lattice $L(A)$ is even and without roots.
It has rank $\vert A\vert -1$ and determinant $\vert A\vert^3$.
The semidirect product $\mathrm{Aut}(A)\ltimes A$ acts isometrically on $L(A)$
in the obvious way. Vectors of norm $4$ in $L(A)$ determine the group $A$ 
uniquely as follows: An arbitrary index of a basis element can be chosen
as the identity $0$
of $A$. A vector $e_0-e_a-e_b+e_c$ yields the identity $a+b=c$
in $A$.

The number of pairs of minimal vectors of norm $4$ in
$L(A)$ is given by the following result:

\begin{prop}\label{propmpLA}
The number of pairs or vectors of norm $4$ in $L(A)$ is given by
$$\vert A\vert
\left(1-\frac{1}{2^c}\right){\vert A\vert/2 \choose 2}+
\frac{\vert A\vert}{2^c}{\left(\vert A\vert-2^c\right)/2\choose 2}$$
where $c$ is the minimal number of generators of the $2$-torsion
subgroup in $A$. Equivalently, $c$ is the largest integer such that 
$A$ contains a subgroup isomorphic to the $c$-dimensional vector space
$\mathbb F_2^c$ over the field $\mathbb F_2$ of two elements.
\end{prop}

\noindent{\bf Proof} We count for each element $a$ of $A$ 
the number $N_a$ 
of solutions of the equation $x+y=a$ with $x,y$ two different elements in $A$.
The total number of pairs of vectors of norm $4$ in $L(A)$
is then given by $\sum_{a\in A}{N_a/2\choose 2}$ since such pairs
are given by $\pm\left(e_{x_1}+e_{y_1}-\left(e_{x_2}+e_{y_2}\right)\right)$
with $\{x_1,y_1\}\not=\{x_2,y_2\}$ such that we have the
equality $x_1+y_1=x_2+y_2$ in $A$.

The kernel of the endomorphism of $A$ defined by $x\longmapsto 2x$
is an $\mathbb F_2$-vector space of dimension $c$. We denote by 
$2A$ its image (of size $\frac{\vert A\vert}{2^c}$) in $A$.
For an element $a\in A\setminus(2A)$ there are $N_a=\vert A\vert$ solutions 
to the equation $x+y=a$ with $x\not=y$ and there are
$\vert A\vert\left(1-\frac{1}{2^c}\right)$ elements in 
$A\setminus (2A)$.

If $a$ is one of the $\frac{\vert A\vert}{2^c}$ elements in $2A$,
there are $2^c$ solutions of $2x=a$ and the equation $x+y=a$
has thus only $N_a=\vert A\vert -2^c$ solutions with $x$ different from
$y$.\hfill$\Box$

\begin{rem}
The algebraic identity
$$\vert A\vert
\left(1-\frac{1}{2^c}\right)\vert A\vert/2+
\frac{\vert A\vert}{2^c}\left(\vert A\vert-2^c\right)/2
={\vert A\vert\choose 2}$$
encodes the fact that $A$ contains ${\vert A\vert\choose 2}$
pairs of distinct elements.
\end{rem}

\begin{thm} \label{thmabelian}
The lattice $L(A)$ associated to an abelian group
having at least $9$ elements is perfect.
\end{thm}

Some lattices $L(A)$ associated to abelian groups $A$ with less than
$9$ elements are perfect.
The lattice $L((\mathbb Z/4\mathbb Z)\oplus
(\mathbb Z/2\mathbb Z))$ is however not perfect (the other
two abelian groups with $8$ elements and the cyclic group
with $7$ elements give rise to perfect lattices, see Sections
\ref{subsectdim7abel} and \ref{subsectdim6abel}).

Given a subset $\mathcal A$ of a finite abelian group $A$,
we define the lattice $L(\mathcal A)$ as the sublattice of 
$L(A)$ generated by all vectors of $L(A)$ involving no elements
of $A\setminus\mathcal A$.

We have the following generalization of Theorem \ref{thmabelian}:

\begin{thm} \label{thmabeliangen}
For a fixed integer $k$ there are only finitely many
isomorphism classes of 
pairs $(\mathcal A\subset A)$ where $A$ is a finite abelian 
group and where $\mathcal A$ is a subset of $A$ with $A\setminus 
\mathcal A$ containing at most $k$ elements such that
$L(\mathcal A)$ is not perfect.
\end{thm}

\subsection{Proofs}

\noindent{\bf Proof of Theorem \ref{thmabelian}}
We establish Theorem \ref{thmabelian} first for cyclic groups.
It holds for $A=\mathbb Z/7\mathbb Z$ and $A=\mathbb Z/8\mathbb Z$
by a direct computation left to the reader (see also Sections 
\ref{subsectdim6abel} and \ref{subsectdim7abel}).
For $N\geq 9$ the $(N-1)$-dimensional lattice $L(\mathbb Z/N\mathbb Z)$
contains the perfect lattice $L_{N-2}$ as a sublattice, 
see Theorem \ref{thmperfectionLdfamily}. By Proposition 
\ref{propfondperf} we need 
to show that minimal vectors of $L(\mathbb Z/N\mathbb Z)$
not orthogonal to $(0,1,2,\dots,N-1)$ span the $(N-1)$-dimensional 
vector space $(1,1,\dots,1)^\perp$.

We consider first the $N-3$ minimal vectors 
\begin{eqnarray*}
v_2&=&e_0+e_1-e_2-e_{N-1},\\
v_3&=&e_0+e_2-e_3-e_{N-1},\\
v_4&=&e_0+e_3-e_4-e_{N-1},\\
&\vdots&\\
v_{N-3}&=&e_0+e_{N-4}-e_{N-3}-e_{N-1},\\
v_{N-2}&=&e_0+e_{N-3}-e_{N-2}-e_{N-1}\\
\end{eqnarray*}
defining the rows of the $(N-3)\times N$ matrix
$$M=\left(\begin{array}{cccccccccccccccccc}
1&1&-1&0&0&0&&&&0&0&0&-1\\
1&0&1&-1&0&0&&&&0&0&0&-1\\
1&0&0&1&-1&0&&&&0&0&0&-1\\
\vdots&&&&\ddots&\ddots&&&&&&&\vdots\\
1&0&0&0&0&0&&&&1&-1&0&-1\\
1&0&0&0&0&0&&&&0&1&-1&-1\\
\end{array}\right)$$
which has obviously rank $N-3$ (column indices are the representatives
$0,1,\dots,N-1$ of $\mathbb Z/N\mathbb Z$).
It is easy to check that $M$ (acting on row-vectors)
has a kernel spanned by the all one vector $(1,1,\dots,1,1)\in\mathbb Z^N$ 
and by the two elements 
\begin{eqnarray*}
w_1&=&(1,0,0,0,0,\dots,0,0,0,1),\\
w_2&=&(1,2,3,\dots,N-2,N-1,0)
\end{eqnarray*}
of $\mathbb Z^N$.
We consider now two additional minimal vectors with signed 
index-sum $N$ given by
\begin{eqnarray*}
v_0&=&e_0+e_1-e_3-e_{N-2},\\
v_1&=&e_1+e_2-e_4-e_{N-1}.
\end{eqnarray*}
Since
$$\left(\begin{array}{cc}
\langle w_1,v_0\rangle&\langle w_1,v_1\rangle\\
\langle w_2,v_0\rangle&\langle w_2,v_1\rangle\\
\end{array}\right)=
\left(\begin{array}{cc}
1&-1\\
1-N&0\\
\end{array}\right)$$
is invertible, the vectors $v_0,\dots,v_{N-1}$ are linearly independent.

In the general case we have to show that linear combinations of
rank $1$ matrices with coefficients $v_av_b,a,b\in A$ for $v=
(v_a)_{a\in A}\in L(A)_{\min}$ have arbitrary off-diagonal coefficients.
Let $(a,b)$ be the index of such an off-diagonal coefficient.
By translation-invariance we can suppose $a=0$.
If $b$ is contained in a cyclic group of order $\geq 7$ we are in the previous
case. We can thus assume that the cardinality of $A$ is divisible 
only by primes $\leq 5$.

If $5$ and either $2$ or $3$ divide the cardinality of $A$, then every
non-zero element of $A$ is contained in a cyclic subgroup of order
at least $10$ and we are done. Otherwise, a non-trivial element of 
$A$ is either contained in a cyclic group of order $25$ (and we are done)
or in $(\mathbb Z/5\mathbb Z)\oplus (\mathbb Z/5\mathbb Z)$
and $L\left((\mathbb Z/5\mathbb Z)\oplus (\mathbb Z/5\mathbb Z)\right)$ 
is perfect by a direct computation (using a Computer-Algebra system).

We are left with the remaining cases where every cyclic subgroup
containing $b$ is of order $2,3,4$ or $6$.
If $b$ is only contained in a cyclic group of order $2$, 
the result follows from perfection of the two groups
$L((\mathbb Z/2\mathbb Z)\oplus(\mathbb Z/2\mathbb Z)\oplus
(\mathbb Z/2\mathbb Z))$ and $L((\mathbb Z/2\mathbb Z)\oplus
(\mathbb Z/8\mathbb Z)$. 
If $b$ is only contained in a cyclic group of order $3$, the result
follows from perfection of 
$L((\mathbb Z/3\mathbb Z)\oplus(\mathbb Z/3\mathbb Z))$.
If $b$ is contained in a cyclic group of order $6$, the result follows 
from perfection of 
$L((\mathbb Z/6\mathbb Z)\oplus(\mathbb Z/2\mathbb Z))$ and 
$L((\mathbb Z/6\mathbb Z)\oplus(\mathbb Z/3\mathbb Z))$.\hfill$\Box$

\begin{prop}\label{propabelcycl} 
If $N$ is large enough then $L(\mathbb Z/N\mathbb Z\setminus
\{a_1,\dots,a_k\})$ is perfect for every subset $\{a_1,\dots,a_k\}$
of $k$ elements in $\mathbb Z/N\mathbb Z$.
\end{prop}

\noindent{\bf Proof} $L(\mathbb Z/N\mathbb Z\setminus
\{a_1,\dots,a_k\})$ contains the lattice $L_{N-k-2}(a_1+1,\dots,a_k+1)$
as a sublattice (we represent elements of $\mathbb Z/N\mathbb Z$ by
natural integers in $\{0,\dots,N-1\}$) and this sublattice is perfect for 
$N\geq \max(9+k,2(k+1)^3+k+1)$ by Theorem \ref{thmperfectfamgen}. 
It is thus enough to show that minimal vectors with signed 
indices summing up to $N$ generate the whole vector-space
$L(\mathbb Z/N\mathbb Z\setminus
\{a_1,\dots,a_k\})\otimes_{\mathbb Z}\mathbb R$.
This can be done (with an effective lower bound on $N$)
as in the proof of Theorem \ref{thmoddperfectfamgen}.
\hfill$\Box$

\begin{prop} \label{propabeliannosmallorder}
There exists an integer $N=N_k$ such that $L(\mathcal A)$ is perfect 
if the finite abelian group $A$ containing $\mathcal A$ 
has an element of prime-order at least $N$ 
and if $A\setminus \mathcal A$ has at most $k$ elements.
\end{prop}

\noindent{\bf Proof} We identify tensor products $v\otimes v$ defined by
elements $v$ in $L(\mathcal A)$ with symmetric matrices whose rows and columns
are indexed by $\mathcal A$. It is enough to show that all such matrices 
with exactly two non-zero diagonal entries and two off-diagonal
non-zero entries defining a symmetric submatrix of the form
$\left(\begin{array}{cc}-1&1\\1&-1\end{array}\right)$
are sums of symmetric matrices associated to 
minimal elements in $L(\mathcal A)$. Up to a translation (of 
$\mathcal A$ and all indices)
we can assume that the first diagonal entry is associated to 
the trivial element $0$ in $A$. The second diagonal element is
then associated to a certain non-zero element $b\in A$ contained in a
cyclic group of order at least $N$ 
and we are done by Proposition \ref{propabelcycl}.
\hfill$\Box$.

\noindent{\bf Proof of Theorem \ref{thmabeliangen}}
As in the proof of Proposition \ref{propabeliannosmallorder}
we want to realize a symmetric matrix corresponding to
$-e_0\otimes e_0+e_0\otimes e_b+e_b\otimes e_0-e_b\otimes e_b$ (up to
a suitable translation), perhaps modulo diagonal matrices.
In particular, we can suppose that $\mathcal A$ contains the trivial 
element $0$.
Proposition \ref{propabelcycl} shows that we can assume that 
every cyclic group containing $b$ is small. 
The group $A$ (if it is huge) has then
a huge number of distinct subgroups.
In particular, we can suppose that it contains a non-trivial 
translate $b+B\not=B$ of a group $B$ containing $a$
with $L(B)$ perfect (this is the case if $B$ has at least $9$
elements by Theorem \ref{thmabelian}). We may now consider the
symmetric matrix $P$ associated to the tensor-product
$$v_1\otimes v_1+v_2\otimes v_2+\dots+v_\alpha\otimes v_\alpha$$
where $\alpha$ is the order of $a$ and where
$$v_i=e_0-e_a+e_{b+(i+1)a}-e_{b+ia}$$
for $i=1,\dots,\alpha$. We have $P_{0,a}=P_{a,0}=\alpha$
and all other non-zero coefficients of $P$ are either diagonal
or have both indices in $b+B$. Coefficients of the last form can be
killed using perfection of $L(B)$. \hfill$\Box$

\begin{rem} Our proof of Theorem \ref{thmabeliangen} can be
unravelled in order to yield effective bounds on the size of $A$.
\end{rem}

\subsection{Examples}

There are no interesting examples in dimension $<6$.

\subsubsection{Dimension $6$}\label{subsectdim6abel}

The $6$-dimensional lattice $L({\mathbb Z/7\mathbb Z})$
associated to the unique group with seven elements has $21$
pairs of minimal elements and is perfect.
A basis is given by the six rows of 
$$A=\left(\begin{array}{ccccccc}
0&1&-1&0&0&-1&1\\
0&1&0&-1&-1&0&1\\
1&1&0&-1&0&-1&0\\
0&1&-1&0&-1&1&0\\
0&1&-1&-1&1&0&0\\
1&0&-1&0&-1&0&1
\end{array}\right)$$
(with columns indexed by the representatives $0,1,2,3,4,5,6$ 
of $\mathbb Z/7\mathbb Z$). Its Gram matrix is the matrix
$$P_6^5=\left(\begin{array}{cccccc}
4&2&2&1&2&2\\
2&4&2&2&1&2\\
2&2&4&0&2&1\\
1&2&0&4&1&2\\
2&1&2&1&4&0\\
2&2&1&2&0&4
\end{array}\right)$$
at page 381 in Chapter XIV of \cite{M}.

\subsubsection{Dimension 7}\label{subsectdim7abel} 
There are $3$ groups with $8$ elements. 

For the cyclic group $\mathbb Z/8\mathbb Z$ we get
$$\left(\frac{8}{2}{4\choose 2}+4{3\choose 2}\right)=36$$
pairs of minimal vectors in the associated $7$-dimensional lattice
$L({\mathbb Z/8\mathbb Z})$
which is perfect and has a basis given by the seven rows of the matrix
$$A=\left(\begin{array}{cccccccc}
0&0&1&-1&0&-1&1&0\\
0&-1&1&0&0&0&1&-1\\
1&0&1&-1&0&0&0&-1\\
0&-1&1&0&1&-1&0&0\\
0&0&-1&-1&0&0&1&1\\
-1&0&0&1&1&0&0&-1\\
1&1&0&0&-1&-1&0&0\\
\end{array}\right)$$
with associated Gram matrix $AA^t$ the matrix 
$$P_7^5=\left(\begin{array}{ccccccc}
4&2&2&2&1&-1&1\\
2&4&2&2&-1&1&-1\\
2&2&4&1&-1&-1&1\\
2&2&1&4&-1&1&-1\\
1&-1&-1&-1&4&-2&0\\
-1&1&-1&1&-2&4&-2\\
1&-1&1&-1&0&-2&4
\end{array}\right)$$
of page 382 in \cite{M}.

The lattice $L((\mathbb Z/4\mathbb Z)\oplus(\mathbb Z/2\mathbb Z))$
(with $38$ pairs of minimal vectors) has perfection-default $2$ and is thus
not perfect.

The lattice $L({(\mathbb Z/2)\oplus(\mathbb Z/2\mathbb Z)\oplus(\mathbb Z/2\mathbb Z)})=
L({\mathbb F_2^3})$ with $42$ pairs of minimal vectors
has a basis given by the seven last rows of the 
table
$$\begin{array}{cccccccc}
000&001&010&011&100&101&110&111\\
\hline
0&0&1&-1&1&-1&0&0\\
-1&1&0&0&1&-1&0&0\\
0&1&1&0&0&-1&-1&0\\
0&0&1&-1&0&0&-1&1\\
-1&0&1&0&1&0&-1&0\\
0&1&0&-1&1&0&-1&0\\
0&0&0&0&1&-1&-1&1
\end{array}$$
(with the first row showing all elements of $\mathbb F_2^3$
corresponding to column-indices). The associated Gram matrix
has only even entries. Dividing it by $2$ we get the matrix
$$P_7^4=\left(\begin{array}{ccccccc}
2&1&1&1&1&1&1\\
1&2&1&0&1&1&1\\
1&1&2&1&1&1&1\\
1&0&1&2&1&1&1\\
1&1&1&1&2&1&1\\
1&1&1&1&1&2&1\\
1&1&1&1&1&1&2
\end{array}\right)$$
(see page 382 in \cite{M}) defining the root lattice $\mathbb D_7$.

\begin{rem} Even parity of all scalar products between minimal vectors
fails for the lattices $L(\mathbb F_2^k)$ with $k\geq 4$.
\end{rem}

\subsubsection{Dimension $8$}

Both $8$-dimensional lattices $L(\mathbb Z/9\mathbb Z)$
and $L((\mathbb Z/3\mathbb Z)\oplus(\mathbb Z/3\mathbb Z))$
have $54$ pairs of minimal vectors and are perfect.
They are non-isomorphic: Every pair of minimal vectors
is orthogonal to exactly 
$15$ pairs of minimal vectors in  $L(\mathbb Z/9\mathbb Z)$ 
and every such pair is orthogonal to exactly 
$9$ pairs of minimal vectors in  $L((\mathbb Z/3\mathbb Z)\oplus
(\mathbb Z/3\mathbb Z))$.
 
\subsection{Examples with one missing element}

The obvious action of $A$ on itself shows that all the lattices
$L(\mathcal A)$ are isomorphic if $\mathcal A$ is obtained by removing
a unique element from $A$. The lattice $L(A\setminus\{0\})$ has
$$\vert A\vert
\left(1-\frac{1}{2^c}\right){\vert A\vert/2-1 \choose 2}+
\frac{\vert A\vert-2^c}{2^c}{\left(\vert A\vert-2^c\right)/2-1\choose 2}
+{\left(\vert A\vert-2^c\right)/2\choose 2}$$
pairs of minimal vectors with norm $4$ 
(where $c$ denotes the dimension of the 
maximal $\mathbb F_2$-vector space occurring as a subgroup in $A$).


\subsection{The root lattice $A_6$}\label{subsectionA6}

Working with the set $\mathcal A=\{001,010,011,100,101,110,111\}$
of all seven non-zero elements in $\mathbb F_2^3$
we get the perfect rescaled root lattice $A_6$ generated by the last seven rows
(with the first row indicating the index-set $\mathcal A$)
of 
$$\begin{array}{ccccccc}
001&010&011&100&101&110&111\\
\hline
0&0&0&1&1&-1&-1\\
0&1&1&0&0&-1&-1\\
-1&0&1&1&0&-1&0\\
-1&1&0&0&1&-1&0\\
-1&1&0&1&0&0&-1\\
-1&0&1&0&1&0&-1
\end{array}$$
Identifying the seven elements of $\mathcal A$ in the obvious way
with the seven points of the Fano plane (projective
plane over $\mathbb F_2)$ we can consider  
pairs of minimal vectors of $L(\mathcal A)$ (i.e. pairs of opposite roots
of $A_6$) as projective lines endowed with marked points (or, dually, 
as points together with incident lines) as follows:
The two coordinates corresponding to coefficients $1$ and the two coordinates 
corresponding to coefficients $-1$
of a minimal vector define projective lines which meet at a point on 
the projective
line associated to the three coordinates corresponding to coefficients $0$.
Up to multiplication by $-1$, this construction is one-to-one
and yields the $21=7\times 3$ pairs of roots of $A_6$.

The Gram matrix associated to the basis of $L(\mathcal A)$ given above
is twice the matrix
$$P_6^7=\left(\begin{array}{cccccc}
2&1&1&1&1&1\\
1&2&1&1&1&1\\
1&1&2&1&1&1\\
1&1&1&2&1&1\\
1&1&1&1&2&1\\
1&1&1&1&1&2\\
\end{array}\right)$$
(which is the Gram matrix with respect to 
the basis $e_0-e_1,\dots,e_0-e_6$ of $A_6$) 
in Chapter XIV of \cite{M}.

\subsection{Two perfect examples of dimension $7$}
Working with $\mathcal A=\{1,\dots,8\}\subset\mathbb Z/9\mathbb Z$,
we get a perfect $7$-dimensional lattice $L(\mathcal A)$ with $30$ pairs 
of minimal vectors. A basis is given by the seven rows of the matrix
$$A=\left(\begin{array}{cccccccc}
1&0&1&0&-1&0&0&-1\\
0&-1&1&1&-1&0&0&0\\
0&-1&1&0&0&0&1&-1\\
1&0&0&1&0&-1&0&-1\\
1&-1&0&0&0&-1&1&0\\
0&1&1&0&0&-1&0&-1\\
0&0&0&1&-1&-1&1&0
\end{array}\right)$$
(with column-indices representing $1,\dots,8\in\mathbb Z/9\mathbb Z$)
with associated Gram matrix $AA^t$ given by the matrix 
$$P_7^{28}=\left(\begin{array}{ccccccc}
4&2&2&2&1&2&1\\
2&4&2&1&1&0&2\\
2&2&4&1&2&1&1\\
2&1&1&4&2&2&2\\
1&1&2&2&4&0&2\\
2&0&1&2&0&4&1\\
1&2&1&2&2&1&4
\end{array}\right)$$ 
(with determinant $2^3\cdot 3^4$) of Chapter XIV in \cite{M}.

The last seven rows of the table
$$\begin{array}{cccccccc}
01&02&10&11&12&20&21&22\\
\hline
1&-1&1&0&-1&0&0&0\\
1&0&1&0&0&-1&-1&0\\
0&0&1&-1&0&0&-1&1\\
1&0&0&-1&-1&0&0&1\\
1&-1&0&0&0&0&-1&1\\
0&0&1&0&-1&-1&0&1\\
0&0&0&-1&1&-1&0&1\\
\end{array}$$
(the first row displays the column indices $\alpha\beta$
with $(\alpha,\beta)\in(\mathbb Z/3\mathbb Z)^2\setminus\{0,0\}$)
define the perfect $7$-dimensional lattice $L\left(
(\mathbb Z/3\mathbb Z)+(\mathbb Z/3\mathbb Z)\setminus\{0,0\}\right)$
with $30$ pairs
of minimal vectors. Its  Gram matrix is
$$P_7^{27}=\left(\begin{array}{ccccccc}
4&2&1&2&2&2&-1\\
2&4&2&1&2&2&1\\
1&2&4&2&2&2&2\\
2&1&2&4&2&2&1\\
2&2&2&2&4&1&1\\
2&2&2&2&1&4&1\\
-1&1&2&1&1&1&4
\end{array}\right)$$
in Chapter XIV of \cite{M}.

\section{The even sublattice construction for abelian groups}
\label{sectevanabelconstr}

Given a finite abelian group $A$ indexing the coordinates of 
$\mathbb Z^A$, we denote by $M(A/(\pm 1))$ the even sublattice 
of $\mathbb Z^{A/(\pm 1)}$ consisting of all elements $v=(v_a)_{a\in A/(\pm a)}$ 
such that $\sum_{a\in A/(\pm 1)}v_a\equiv 0\pmod 2$ (this ensures
evenness of $M(A/(\pm 1))$ and such that $\sum_{a\in A/(\pm 1)}v_aa=0\in A$ where
$A/(\pm 1)$ denotes (somewhat abusively) 
a set of representatives of $A$ under the involutive automorphism 
$a\longmapsto -a$. The lattice $M(A/(\pm 1))$ is
without roots. It has rank $\vert A/(\pm 1)\vert$ and determinant
$4\vert A\vert^2$.
Vectors of norm $4$ in $M(A/(\pm 1))$ are of the form $\pm 2e_a$ if $2a=0$
in $A$ for $a\in A/(\pm 1)$ or of the form 
$\pm e_{a_1}\pm e_{a_2}\pm e_{a_3}\pm e_{a_4}$ if 
$\pm {a_1}\pm {a_2}\pm {a_3}\pm {a_4}=0$ in $A$ for four distinct
elements $a_1,\dots,a_4$ of $A/(\pm 1)$ with $\pm$ denoting suitable 
choices of signs. The subgroup of all elements of order at most $2$ acts
by isometries on $M(A/(\pm 1))$ and the group $A$ can be recovered 
(up to isometries) from the set of minimal vectors of norm $4$
in $M(A/(\pm 1))$.

\begin{thm}\label{thmevsublattcycl} The lattice $M((\mathbb Z/N\mathbb Z)/(\pm 1))$
associated to a cyclic group of order $\geq 15$ is perfect.
\end{thm}

Theorem \ref{thmevsublattcycl} can probably be generalized to arbitrary
finite abelian groups which are sufficiently large. It 
should have a further generalization 
obtained by removing $k$ elements from $A/(\pm 1)$.

\noindent{\bf Proof of Theorem \ref{thmevsublattcycl}}
We consider first a cyclic group $A=\mathbb Z/N\mathbb Z$ of even order
$N=2m$. Representatives of $A/(\pm 1)$ are
$\{0,1,\dots,m\}$. For $N\geq 16$, the lattice $M(A/(\pm 1))$ contains 
the perfect sublattice $M_{m-1}=M(A/(\pm 1))\cap (0,1,2,\dots,m)^\perp$, 
see Theorem \ref{thmevensublconstr}.

We set $v_i=-e_i+e_{i+1}+e_{m-1}+e_m$ for $i=0,\dots,m-3$.
The minimal elements $v_0,\dots,v_{m-3}$ together with $2e_0=(e_0+e_1+e_{m-1}+e_m)-(-e_0+e_1+e_{m-1}+e_m),2e_m$
and $e_1+e_2+e_{m-2}+e_{m-1}$ are linearly independent. Since the corresponding
signed index-sum $-i+i+(m-1)+m=2m$ (respectively $0+1+(m-1)+m=2m$ and $2m$)
is non-zero they are not orthogonal to $(0,1,2,\dots,m)$.
Perfection of $M((\mathbb Z/N\mathbb Z)/(\pm 1))$ for even $N\geq 16$ 
follows now from  Proposition \ref{propfondperf}.

For a cyclic group $N=2m+1$ of odd order $2m+1$ we proceed as follows:
The $m-3$ linearly independent minimal 
elements $-e_i+e_{i+2}+e_{m-1}+e_{m},i=0,\dots,m-4$
can be completed to a base by adjoining the following four elements
\begin{eqnarray*}
u_1&=&(1,0,1,0,1,0,1,0,\dots,0,0)\\
u_2&=&(1,1,1,1,\dots,1,1,0,0)\\
u_3&=&(0,1,2,3,4,5,\dots,m-3,m-2,-1,-1)\\
u_4&=&(0,0,0,0,0,\dots,0,0,0,1,-1)
\end{eqnarray*}
($u_1$ has alternating coefficients $0,1$ except for the last
two coefficients which are both zero) which are orthogonal
to $-e_i+e_{i+2}+e_{m-1}+e_{m}$ for $i\in\{0, \dots,m-4\}$.
We consider now four minimal vectors given by
\begin{eqnarray*}
w_1&=&e_0+e_2+e_{m-1}+e_m\\
w_2&=&e_0+e_3+e_{m-2}+e_m\\
w_3&=&e_0+e_4+e_{m-2}+e_{m-1}\\
w_4&=&e_1+e_3+e_{m-2}+e_{m-1}
\end{eqnarray*}
The matrix $S$ of scalar products $S_{i,j}=\langle w_i,u_j\rangle$ equals 
$$\left(\begin{array}{cccc}
2&2&0&0\\
1+\epsilon&3&m&-1\\
2+\epsilon&3&m+1&1\\
3+\epsilon&3&m+1&1\end{array}\right)$$
where $\epsilon=0$ if $m$ is odd and $\epsilon=1$ if $m$ is even.
The matrix $S$ has non-zero determinant $8m+4$ which ends the proof
by Proposition \ref{propfondperf}.
\hfill$\Box$

\subsection{A non-cyclic example giving $E_8$}\label{subsectionE8}

All elements of the additive group $\mathbb F_2^3$ are their own inverses
and $M((\mathbb F_2)^3)/(\pm 1))$ is obtained from 
the lattice $L(\mathbb F_2^3)$ by considering 
$L(\mathbb F_2^3)+(2\mathbb Z)^{\mathbb F_2^3}$. The resulting lattice is
the (rescaled) exceptional root-lattice $E_8$
with basis the last eight rows of
$$\begin{array}{cccccccc}
000&001&010&011&100&101&110&111\\
\hline
2&0&0&0&0&0&0&0\\
-1&-1&-1&-1&0&0&0&0\\
0&2&0&0&0&0&0&0\\
0&-1&1&0&-1&0&0&-1\\
0&0&0&0&1&-1&1&1\\
0&0&0&0&0&2&0&0\\
0&0&0&0&-1&-1&-1&1\\
0&0&-1&1&0&0&-1&-1\\
\end{array}$$
having twice the Dynkin matrix
$$\left(\begin{array}{cccccccc}
2&-1&0&0&0&0&0&0\\
-1&2&-1&0&0&0&0&0\\
0&-1&2&-1&0&0&0&0\\
0&0&-1&2&-1&0&0&0\\
0&0&0&-1&2&-1&0&-1\\
0&0&0&0&-1&2&-1&0\\
0&0&0&0&0&-1&2&0\\
0&0&0&0&-1&0&0&2\\
\end{array}\right)$$
of $E_8$ as its Gram matrix . 

\subsection{Removing an element}\label{subsectionE7}
(One can in fact remove an arbitrary element from $\mathbb F_2^3$.)
The even lattice associated to all $7$
non-zero elements $\mathbb F_2^3\setminus\{0\}$ of $\mathbb F_2^3$ is
the lattice generated by the seven vectors
$$\begin{array}{ccccccc}
001&010&011&100&101&110&111\\
\hline
0&0&0&1&1&1&-1\\
0&0&0&2&0&0&0\\
1&0&1&1&0&1&0\\
0&0&0&1&1&-1&-1\\
0&0&0&1&1&1&1\\
0&1&1&1&1&0&0\\
0&-1&1&1&1&0&0.\\
\end{array}$$
The associated Gram matrix is twice the matrix
$$P_7^1=
\left(\begin{array}{rrrrrrr}
2&1&1&1&1&1&1\\
1&2&1&1&1&1&1\\
1&1&2&0&1&1&1\\
1&1&0&2&0&1&1\\
1&1&1&0&2&1&1\\
1&1&1&1&1&2&1\\
1&1&1&1&1&1&2
\end{array}\right),$$
see \cite{M}, page 382, defining the exceptional root lattice 
$E_7$.

\section{A construction 
with minimum $3$ using $\mathbb F_2^c$}\label{sectmin3}

Given a finite-dimensional vector space $\mathbb F_2^c$ of dimension $c$
over the field $\mathbb F_2$ of two elements, the lattice $T(\mathbb F^c)
$ is the integral
sublattice of $\mathbb Z^{\mathbb F^c\setminus\{0\}}$ consisting of all vectors
$v=(v_a)_{a\in\mathbb F_2^c\setminus\{0\}}$ such that 
$\sum_{a\in \mathbb F_2^c\setminus\{0\}}
v_aa=0$ in $\mathbb F_2^c$. Minimal vectors have norm $3$ (except 
in the trivial case $c=1$) and are given by
$\epsilon_1e_a+\epsilon_2e_b+\epsilon_3e_c$
with $\epsilon_1,\epsilon_2,\epsilon_3\in\{\pm 1\}$
and with $a,b,c=a+b\in \mathbb F_2^c$ defining a projective line
of the $(c-1)$-dimensional projective space over $\mathbb F_2$.

\begin{thm} \label{thmconstrT}
The lattice $T(\mathbb F_2^c)$ has no roots, determinant $4^c$ 
and $\frac{4}{3}{2^c-1\choose 2}$
pairs of vectors of norm $3$. It is perfect for $c\geq 3$.
\end{thm}

\noindent{\bf Proof}
The lattice $T(\mathbb F_2^c)$ is the kernel of the augmentation-map. 
It is thus of index $2^c$ in $\mathbb Z^{2^c-1}$
and has determinant $(2^c)^2=4^c$. There are 
$\frac{(2^c-1)(2^c-2)}{3\cdot 2}=\frac{1}{3}{2^c-1\choose 2}$ 
projective lines in $\mathbb F_2^c\setminus\{0\}$ and every 
projective line determines $4$ pairs of minimal 
vectors. 

In order to prove perfection, 
we consider a symmetric matrix $S$ with $2^c-1$ rows and columns indexed 
by all non-zero elements of $\mathbb F_2^c$.
A non-zero diagonal coefficient $c_{a,b}$ of $S$ can be eliminated by 
subtracting 
$$\frac{c_{a,b}}{4}\left(-v_{+++}v_{+++}^t-v_{+-+}v_{+-+}^t+
v_{++-}v_{++-}^t+v_{+--}v_{+--}^t\right)$$
from $S$ where $v_{\epsilon_1,\epsilon_2\epsilon_3}=\epsilon_1e_a
+\epsilon_2e_b+\epsilon_3e_c$ with $c=a+b\in \mathbb F_2^c$.

The orthogonal projector 
$$\frac{1}{4}\left(v_{+++}v_{+++}^t+v_{+-+}v_{+-+}^t+
v_{++-}v_{++-}^t+v_{+--}v_{+--}^t\right)$$
has only three non-zero coefficients on the diagonal corresponding
to rows (and columns) indexed by $a,b$ and $c=a+b$. It is thus 
associated to the diagonal coefficient of a projective line 
over $\mathbb F_2$.
The matrix $A$ defined by the last seven rows (with the first row indicating
the seven points of the projective plane over $\mathbb F_2$) of 
$$\begin{array}{ccccccc}
001&010&011&100&101&110&111\\
\hline
1&1&1&0&0&0&0\\
1&0&0&1&1&0&0\\
1&0&0&0&0&1&1\\
0&1&0&1&0&1&0\\
0&1&0&0&1&0&1\\
0&0&1&1&0&0&1\\
0&0&1&0&1&1&0\\
\end{array}$$
has determinant $-24$ and is thus invertible.
This shows that we can get rid of diagonal coefficients
using the ``diagonal'' projectors onto projective lines
by embedding them into projective planes. More precisely,
given a point $a\in\mathbb F_c$ of a projective
plane $\Pi$, the projector 
$$\frac{1}{6}\left(2\sum_{a\in L\subset \Pi}P_L-\sum_{a\not\in L\subset\Pi}P_L\right)$$
is the diagonal projector onto the diagonal element indexed by $a$
where $P_L$ is the projector 
$$\frac{1}{4}\left(v_{+++}v_{+++}^t+v_{+-+}v_{+-+}^t+
v_{++-}v_{++-}^t+v_{+--}v_{+--}^t\right)$$
(with $v_{\epsilon_1\epsilon_2\epsilon_3}$ as above)
associated to projective line $\{a,b,c=a+b\}\subset \Pi$.
\hfill$\Box$

\begin{rem} (i) No elements (except subsets leaving the non-zero elements
of a subgroup containing at least $8$ elements) 
can be removed from the set $\mathbb F_2^c\setminus\{0\}$
in the construction of $T(\mathbb F_2^c)$
without destroying perfection of the associated lattice.

The construction cannot be adapted to other finite abelian groups
(with $\mathbb F_2^c\setminus\{0\}$ replaced by representatives
of all non-zero orbits of a finite abelian group 
$A$ under the automorphism $x\longmapsto -x$) without losing 
perfection.
\end{rem}

\subsection{Digression: The equiangular system of 
the perfect lattice $T(\mathbb F_2^3)$ and the Schl\"afli graph}

The $7$-dimensional 
perfect lattice $T(\mathbb F_2^3)$ with $28$ pairs of minimal vectors 
has a basis given by the last seven rows of 
$$\begin{array}{ccccccc}
001&010&011&100&101&110&111\\
\hline
0&0&1&0&1&-1&0\\
0&1&0&0&1&0&-1\\
1&0&0&0&0&-1&1\\
1&0&0&-1&1&0&0\\
0&0&1&-1&0&0&1\\
0&0&1&1&0&0&-1\\
-1&0&0&1&1&0&0
\end{array}$$
with Gram matrix
$$P_7^2=\left(\begin{array}{ccccccc}
3&1&1&1&1&1&1\\
1&3&-1&1&-1&1&1\\
1&-1&3&1&1&-1&-1\\
1&1&1&3&1&-1&-1\\
1&-1&1&1&3&-1&-1\\
1&1&-1&-1&-1&3&1\\
1&1&-1&-1&-1&1&3
\end{array}\right),$$
see page 382 of \cite{M}. Up to rescaling, this is the dual lattice 
$E_7^*$ of the root lattice $E_7$.

Its $28$ pairs of minimal vectors define a system of $28$ equiangular
lines (meeting two-by-two in a common angle given by 
$\mathrm{arccos}(1/3)$) in $\mathbb R^7$.
Supports of minimal vectors define projective lines in the Fano plane 
(projective plane over $\mathbb F^2$).
The automorphism group of $T(\mathbb F_2^3)$ acts transitively
on the set of minimal vectors. Fixing a first minimal vector, 
say $w=(1,1,1,0,0,0,0)$, we chose representatives $v_1,\dots,v_{27}$
of the $27$ minimal pairs different from $\pm w$ such that 
$\langle w,v_i\rangle=1$.
We encode the angles between $v_1,\dots,v_{27}$ by a graph
$\Gamma$ with vertices $v_1,\dots,v_{27}$ and edges $v_i,v_j$ if 
$\langle v_i,v_j\rangle=-1$.
The characteristic polynomial of the 
adjacency matrix $A$ of $\Gamma$ (with diagonal zero and 
off-diagonal coefficients $\frac{1-\langle v_i,v_j\rangle}{2}$)
is given by 
$$(t-10)(t-1)^{20}(t+5)^6$$
and the graph $\Gamma$ is thus a strongly regular graph on $27$ 
vertices with parameters $(v,k,\lambda,\mu)=(27,10,1,5)$.

Otherwise stated, the graph $\Gamma$ has $v=27$ vertices.
It is of degree $k=10$ and diameter
$2$ such that two adjacent vertices in $\Gamma$ have always $\lambda=1$ common
neighbours and two non-adjacent vertices of $\Gamma$ have always $\mu=5$
common neighbours. $\lambda=1$ is equivalent to the fact that 
every edge of $\Gamma$ is contained in a unique triangle (complete
graph on $3$ vertices) of $\Gamma$.

Such a graph is unique and it (or sometimes its complement)
is called the Schl\"afli graph.

\begin{rem}\label{remevensublatticeT3} 
(i) The even sublattice of the lattice $T(\mathbb F_3^2)$
is (up to rescaling) the root lattice $E_7$ consisting of all vectors
of the lattice $M(\mathbb F_2^3)$ (see Section \ref{subsectionE7})
not involving the basis vector
$e_0$ associated to the identity $0$ of the additive group 
$\mathbb F_2^3$. Its $63$ pairs of minimal vectors can be described as 
follows: Every line $\{i,j,k\}$ gives rise to $2^3=8$ pairs of minimal
vectors by considering a vector with zero coordinates corresponding to
$i,j,k$ and with coordinates $\pm 1$ associated to points not in 
$\{i,j,k\}$. This gives $7\cdot 8=56$ pairs of minimal vectors (of norm $4$).
Seven additional pairs are given by $\pm 2e_i$ and are associated
to the seven points of the projective plane.

(ii) Restricting to vectors with zero coordinate-sum
of the even sublattice of $T(\mathbb F_2^3)$, 
we get the rescaled root lattice $A_6$ of Section \ref{subsectionA6}.
\end{rem}

\section{Generalizations}

All lattices constructed in this paper are of the form 
$\Lambda=\ker(\varphi(\mathbb Z^{d+a}))$ for a surjective morphism $\varphi$ from 
$\mathbb Z^{d+a}$ onto an abelian group $\mathbb Z^a\oplus A$ with
$A$ finite. A suitable choice of $\varphi$ ensures nice combinatorial
properties of small elements in $\Lambda$. Up to this 
point (except in Section \ref{sectmin3}), 
we have  worked with even lattices containing no roots
and we have used properties of $\varphi$ for proving perfection 
of the set $\Lambda_4$ of minimal vectors in $\Lambda$.
It is of course tempting to consider $\varphi$ such that 
the norm $\lambda_1^2+\dots+\lambda_{d+a}^2$ of every non-zero element 
$(\lambda_1,\dots,\lambda_{d+a})\in \Lambda$ is at least equal to 
some larger integer $m>4$. Sidon sets provide examples leading to 
minimum $6$ (but do not ensure perfection) as follows:
A \emph{Sidon set} in an additive group $A$ is a subset $\mathcal S$
such that $x_1+y_1=x_2+y_2$ implies $\{x_1,y_1\}=\{x_2,y_2\}$ 
as multisets for $x_1,y_1,x_2,y_2\in\mathcal S$.
The sublattice of all elements in $\mathbb Z^{\mathcal S}$
with zero coefficient-sum $\sum_{x\in\mathcal S} \lambda_x=0$ such that 
$\sum_{x\in\mathcal S} \lambda_xx=0\in A$ is then even and without roots
or vectors of norm $4$. More generally, one might consider 
subsets $\mathcal S$ which have the $m$-lattice Sidon property:
every non-zero vector in the lattice of all elements in $\mathbb Z^{\mathcal S}$
with zero coefficient-sum $\sum_{x\in\mathcal S} \lambda_x=0$ such that 
$\sum_{x\in\mathcal S} \lambda_xx=0\in A$ has (squared Euclidean)
norm at least $2(m+1)$.
As a variation, one can drop the requirement
$\sum_{x\in\mathcal S} \lambda_x=0$ by replacing it with 
the evenness condition $\sum_{x\in\mathcal S} \lambda_x\equiv 0\pmod 2$
or by dropping it without any other requirement altogether
(this puts of course an even stronger constraint on $\mathcal S$).

\subsection{Craig lattices}

Given a finite field $\mathbb F_q$ with $q=p^e$ a prime power and an
integer $k$, we can consider the lattice $C_{q-1,k}$ defined by all
vectors of $\mathbb Z^{\mathbb F_q}$ with zero coefficient sum
 $\sum_{x\in\mathbb F_q}\lambda_x=0$
and such that $\sum_{x\in\mathbb F_q}\lambda_x x^i=0\in \mathbb F_q$
for $i=1,\dots,k$ (equality holds of course also for $i=0$).
For $q=p$ a prime number, the lattice $C_{q-1,k}$ is a Craig lattice.
The lattice $C_{q-1,k}$ is even and has determinant $q^{2k+1}$.

\begin{prop}\label{propminCraig} 
The lattice $C_{q-1,k}$ has minimum $\geq 2k+2$ if $k$ is smaller than the 
characteristic $p$ of $\mathbb F_q$. 
\end{prop}

\noindent{\bf Proof} Symmetric power-sums of degree up to $p-1$
define elementary symmetric polynomials of degree up to $p-1$.
A minimal vector with strictly positive coefficients of indices 
$a_1,\dots,a_l$ (with indices repeated according to the value 
of the associated integral coefficient) and strictly negative coefficients of 
indices $b_1,\dots,b_l$ gives rise to two polynomials
$\prod_{i=1}^l(x-a_i)$ and $\prod_{i=1}^l(x-b_i)$. 
Since symmetric power-sums of degree up to $p-1$
define elementary symmetric polynomials of degree up to $p-1$
this implies either of $l>k$ or $k\geq p$.\hfill$\Box$

\begin{prop} For $k$ smaller than the characteristic $p$ of $\mathbb F_q$,
the number of pairs of elements of norm $2(k+1)$ in $C_{q-1,k}$
is given by
\begin{eqnarray}\label{formulaminvectCraig}
\sum_{(a_1,\dots,a_k)\in \mathbb F_q^k}{N(a_1,\dots,a_k)\choose 2}
\end{eqnarray}
where $N(a_1,\dots,a_k)\leq q$ is the number of constants $a_0\in \mathbb F_q$
such that the polynomial 
$x^{k+1}+\sum_{i=0}^ka_ix^{i}$ has exactly $k+1$ distinct roots in 
$\mathbb F_q$. 
\end{prop}

\noindent{\bf Proof} $N(a_1,\dots,a_k)$ is also the number 
of subsets $\{x_1,\dots,x_{k+1}\}$ of $k+1$ distinct elements in 
$\mathbb F_q$ such that $\sum_{i=1}^{k+1}x_i^j=b_j$
with $b_1,\dots,b_k$ the power-sums corresponding to the elementary
symmetric functions $a_k,a_{k-1},\dots,a_1$. Such subsets are disjoint 
and pairs of two such subsets define indices of coefficients $1$ and $-1$
in minimal vectors.\hfill$\Box$

\begin{cor} The lattice $C_{q-1,k}$ (for $k$ smaller than the characteristic
$p$ of $\mathbb F_q$) has at least 
$$q^k{\frac{1}{q^k}{q\choose k+1}\choose 2}$$
pairs of vectors of norm $2(k+1)$.

In particular, for a fixed value of $k$, the lattice $C_{q-1,k}$
has asymptotically at least $\frac{q^{k+2}}{2((k+1)!)^2}$
pairs of minimal vectors of norm $2(k+1)$.
\end{cor}

\noindent {\bf Proof} Since every subset of $k+1$ elements in $\mathbb F_q$
contributes $1$ to exactly one of the numbers $N(a_1,\dots,a_k)$ we have
$${q\choose k+1}=\sum_{(a_1,\dots,a_k)\in \mathbb F_q^k}N(a_1,\dots,a_k).$$
Convexity properties of the polynomial ${x\choose 2}=\frac{x(x-1)}{2}$
imply that (\ref{formulaminvectCraig})
is minimal if all $q^k$ numbers $N(a_1,\dots,a_k)$ are equal.
\hfill$\Box$

\begin{thm} For $k=2$, the number of minimal pairs in $C_{q-1,2}$ is given by
$$\frac{1}{72}q(q-1)(q^2-10q+33)$$
for $q$ a prime-power congruent to $1$ modulo $6$ and by 
$$\frac{1}{72}q(q-1)(q-5)^2$$
for $q$ a prime power congruent to $5$ modulo $6$. 
\end{thm}

\noindent{\bf Sketch of Proof} We have to evaluate Formula
(\ref{formulaminvectCraig}) for $k=2$. Substituting $x$ 
with $x-\frac{a_2}{3}$ we get $N(a_1,a_2)=N(a_1-\frac{1}{3}a_2^2,0)$.
Formula (\ref{formulaminvectCraig}) for $k=2$ is thus given by
$$q\sum_{a\in \mathbb F_q}{N(a,0)\choose 2}$$
if $q$ is not a power of $3$.
Since $N(a,0)$ depends only on the value $\left(\frac{a}{q}\right)$
of the quadratic character extending the Jacobi symbol, we have
to compute $N(a,0)$ for $a=0,1$ and for a non-square of $\mathbb F_q$.
These computations boil down to classical properties of binary quadratic forms
over finite fields. (One can alternatively use a result of Stickelberger,
as observed by the reviewer.)\hfill$\Box$

\begin{rem} A close relative of the lattice $C_{q-1,2}$ is the lattice 
associated to the Sidon set $\{(x,x^{-1})\}_{x\in\mathbb F_q^*}
\subset \mathbb F_q^2$
for $\mathbb F_q$ a finite field of odd characteristic.
It is of dimension $q-2$, has minimum $6$ (except for a few small
values of $q$) and consist of all elements 
$(\lambda_x)_{x\in\mathbb F_q^*}\in\mathbb Z^{\mathbb F_q^*}$ 
(integral vectors with indices in $\mathbb F_q^*$)
such that $\sum_{x\in\mathbb F_q^*}\lambda_x=0$ and 
$\sum_{x\in\mathbb F_q^*}\lambda_xx=\sum_{x\in\mathbb F_q^*}\lambda_xx^{-1}=0
\in\mathbb F_q$.
\end{rem}

For $k=3$, let $c_q$ be such that the number of pairs of minimal vectors 
(of norm $8$) in $C_{q-1,3}$ is given by 
$$\frac{1}{1152}q(q-1)(q^3-21q^2+171q-c_q).$$
Writing $c_k$ as
$$c_k=483+36\left(\frac{-1}{q}\right)+64\left(\frac{-3}{q}\right)+\delta_q,$$
we have the following result due to  Noam D. Elkies, see \cite{MOq}
(a preliminary draft of the present paper
proposed the values corresponding to $\delta_q=0$ conjecturally):

\begin{thm} If $q$ is a prime $\leq 5$, then 
$\delta_q=0$ if $\left(\frac{-2}{q}\right)=-1$
(yielding the values
$$c_q=\left\lbrace\begin{array}{ll}
455\quad&\hbox{ if }q\equiv 5\pmod{24},\\
511\quad&\hbox{ if }q\equiv 7\pmod{24},\\
583\quad&\hbox{ if }q\equiv 13\pmod{24},\\
383\quad&\hbox{ if }q\equiv 23\pmod{24}\end{array}\right.$$
for $c_q$ in these cases) and
$$\delta_q=24(m^2-2n^2)+192+72\left(\frac{-1}{q}\right)$$
where $m$ and $n$ are the unique natural integers such that $q=m^2+2n^2$
otherwise (i.e. for $q\geq 11$ a prime such that $\left(\frac{-2}{q}\right)=1$).
\end{thm}
See \cite{MOq} for the fairly sophisticated proof.



\noindent{\bf Acknowledgements} I thank Philippe Eyssidieux for an
interesting discussion, Jacques Martinet and an anonymous referee
for a careful reading and many helpful remarks.


\noindent Roland BACHER, 

\noindent Univ. Grenoble Alpes, Institut Fourier, 

\noindent F-38000 Grenoble, France.
\vskip0.5cm
\noindent e-mail: Roland.Bacher@ujf-grenoble.fr

\end{document}